\documentclass[11pt]{elsarticle}
\usepackage{graphicx,subfigure}
\usepackage{amsfonts}
\usepackage{mathrsfs}
\usepackage{amssymb}
\usepackage{chngpage}
\usepackage{amsmath}
\usepackage{graphicx}
\usepackage{amsmath,amssymb}
\usepackage[mathscr]{eucal}
\usepackage{xcolor}
\usepackage{lineno,hyperref}

\newtheorem{theorem}{Theorem}

\newtheorem{corollary}[theorem]{Corollary}
\newtheorem{lemma}[theorem]{Lemma}
\newtheorem{remark}[theorem]{Remark}

\numberwithin{equation}{section}
\numberwithin{theorem}{section}

\def\({\left(}
\def\){\right)}

\makeatletter
\def\ps@pprintTitle{%
  \let\@oddhead\@empty
  \let\@evenhead\@empty
  \let\@oddfoot\@empty
  \let\@evenfoot\@oddfoot
}
\makeatother

\begin{document}

\begin{frontmatter}

\title{Asymptotic behavior of a semilinear problem in heat conduction with long time memory and non-local diffusion}


\author[mymainaddress]{Jiaohui Xu}
\ead{jiaxu1@alum.us.es}

\author[mymainaddress]{Tom\'{a}s Caraballo\corref{mycorrespondingauthor}}
\cortext[mycorrespondingauthor]{Corresponding author}
\ead{caraball@us.es}

\author[mysecondaryaddress]{Jos\'e Valero}
\ead{jvalero@umh.es}

\address[mymainaddress]{
Departamento de Ecuaciones Diferenciales y An\'alisis Num\'erico, Facultad de Matem\'aticas, Universidad de Sevilla,  C/Tarfia s/n, 41012-Sevilla, Spain.}
\address[mysecondaryaddress]{Centro de Investigaci\'on Operativa, Universidad Miguel Hern\'andez, Avda. de la Universidad, s/n, 03202-Elche, Spain.}

\begin{abstract}
In this paper,  the asymptotic behavior of a semilinear heat equation with long time memory and non-local diffusion is analyzed in the usual set-up for dynamical systems generated by differential equations with delay terms. This approach is different from the previous published literature on the long time behavior of heat equations with memory which is carried out by the Dafermos transformation. As a consequence, the  obtained results provide complete information about the attracting sets for the original problem, instead of the transformed one. In particular, the proved results also generalize and complete previous literature in the local case.
\end{abstract}


\end{frontmatter}


\section{Introduction}


The main objective of this paper is to analyze the asymptotic behavior of a
semilinear heat equation with long time memory and non-local diffusion,
which is an interesting situation with important applications in the real
world.

On the one hand, the effects that memory terms (or the past history of a
phenomenon) produce on the evolution of a dynamical system is obvious, since
it is sensible to think that the evolution of any system depends not only on
the current state but on its whole history (see, for instance, \cite{CA, FZ,
C, CR, CM, G, L2} and the references therein). On the other hand, many
problems are better described by considering non-local terms, which created
a great interest in the modeling of various real applications (see \cite{C1,
C2, C3, C} and the references therein).

Motivated by some physical problems from thermal memory or materials with
memory, one can find a significant literature devoted to the analysis of
partial differential equations with long time memory. For example, the
authors introduced in \cite{C} a semilinear partial differential equation to
model the heat flow in a rigid, isotropic, homogeneous heat conductor with
linear memory, which is given by 
\begin{equation}
\begin{cases}
c_{0}\partial _{t}u-k_{0}\Delta u-\displaystyle\int_{-\infty
}^{t}k(t-s)\Delta u(s)ds+f(u)=h, \\ 
u(x,t)=0, \\[0.8ex] 
u(x,\tau +t)=u_{0}(x,t), \\ 
\end{cases}%
\begin{aligned} &\mbox{in}~~\Omega\times (\tau,+\infty),\\
&\mbox{on}~\partial\Omega\times(\tau,+\infty),\\[0.8ex]
&\mbox{in}~~\Omega\times(-\infty,0],\\ \end{aligned}    \label{eq1-1}
\end{equation}%
where $\Omega \subset \mathbb{R}^{N}$ is a bounded domain with regular
boundary, $u:\Omega \times \mathbb{R}\rightarrow \mathbb{R}$ is the
temperature field, $k:\mathbb{R}^{+}\rightarrow \mathbb{R}$ is the heat flux
memory kernel, $\mathbb{R}^{+}$ denotes the interval $(0,+\infty )$, $c_{0}$
and $k_{0}$ denote the specific heat and the instantaneous conductivity,
respectively. To solve \eqref{eq1-1} successfully, the authors considered
this problem as a non-delay one by making the past history of $u$ from $%
-\infty $ to $0^{-}$ be part of the forcing term given by the causal
function $g$, which is defined by 
\begin{equation*}
g(x,t)=h(x,t)+\int_{-\infty }^{\tau }k(t-s)\Delta u_{0}(x,s)ds,\qquad x\in
\Omega ,~~t\geq \tau .
\end{equation*}%
In this way, \eqref{eq1-1} becomes an initial value problem without delay or
memory, 
\begin{equation}
\begin{cases}
c_{0}\partial _{t}u-k_{0}\Delta u-\displaystyle\int_{\tau }^{t}k(t-s)\Delta
u(s)ds+f(u)=g, \\ 
u(x,t)=0, \\ 
u(x,\tau )=u_{0}(x,0), \\ 
\end{cases}%
\begin{aligned} &\mbox{in}~~\Omega\times (\tau,+\infty),\\
&\mbox{on}~\partial\Omega\times (\tau,+\infty),\\ &\mbox{in}~~\Omega.\\
\end{aligned}    \label{eq1-1bis}
\end{equation}%
However, this problem does not generate a dynamical system in an appropriate
phase space, since the equation in \eqref{eq1-1bis} depends on the past
history and we are just fixing an initial value at time $\tau $.

Therefore, two alternatives are possible. The first one is based on the idea
introduced by Dafermos \cite{D}, for linear viscoelasticity, in the 70's.
Let us define the new variables, 
\begin{equation*}
u^{t}(x,s)=u(x,t-s),\qquad s\geq 0,~~t\geq \tau ,
\end{equation*}%
\begin{equation}
\eta ^{t}(x,s)=\int_{0}^{s}u^{t}(x,r)dr=\int_{t-s}^{t}u(x,r)dr,\qquad s\geq
0,~~t\geq \tau .  \label{eq3.2a}
\end{equation}%
Besides, assuming $k(\infty )=0$, a change of variable and a formal
integration by parts imply 
\begin{equation*}
\int_{-\infty }^{t}k(t-s)\Delta u(s)ds=-\int_{0}^{\infty }k^{\prime
}(s)\Delta \eta ^{t}(s)ds.
\end{equation*}%
Setting 
\begin{equation*}
\mu (s)=-k^{\prime }(s),
\end{equation*}%
the original equation \eqref{eq1-1bis} turns into the following autonomous
system without delay, 
\begin{equation}
\begin{cases}
c_{0}\frac{\partial u}{\partial t}-k_{0}\Delta u-\displaystyle%
\int_{0}^{\infty }\mu (s)\Delta \eta ^{t}(s)ds+f(u)=g, \\ 
\eta _{t}^{t}(s)=-\eta _{s}^{t}(s)+u(t), \\ 
u(x,t)=\eta ^{t}(x,s)=0, \\ 
u(x,\tau )=u_{0}(0), \\ 
\eta ^{\tau }(x,s)=\eta _{0}(s),%
\end{cases}%
\begin{aligned} &\mbox{in}~~\Omega\times (\tau,\infty),\\
&\mbox{in}~~\Omega\times (\tau,\infty)\times\mathbb{R}^+,\\
&\mbox{on}~~\partial \Omega\times \mathbb{R}\times\mathbb{R}^+,\\
&\mbox{in}~~ \Omega,\\ &\mbox{in}~~\Omega\times\mathbb{R}^+, \end{aligned} 
   \label{eq3.3a}
\end{equation}%
where, $\eta _{s}^{t}$ denotes the distributional derivative of $\eta
^{t}(s) $ with respect to the internal variable $s$. It follows from the
definition of $\eta ^{t}(x,s)$ (see (\ref{eq3.2a})) that 
\begin{equation}
\eta _{0}(s)=\int_{\tau -s}^{\tau }u(r)dr=\int_{\tau -s}^{\tau }u_{0}(r-\tau
)dr=\int_{-s}^{0}u_{0}(r)dr,  \label{eq2.5a}
\end{equation}%
which is the initial integrated past history of $u$ with vanishing boundary.
Consequently, any solution to \eqref{eq1-1bis} is a solution to %
\eqref{eq3.3a} for the corresponding initial values $(u_{0}(0),\eta _{0})$
given by \eqref{eq2.5a}. It is worth emphasizing that problem \eqref{eq3.3a}
can be solved for arbitrary initial values $(u_{0},\eta _{0})$ in a proper
phase space $L^{2}(\Omega )\times L_{\mu }^{2}(\mathbb{R}^{+};H_{0}^{1}(%
\Omega ))$ (see Section~\ref{s2} for more details), i.e., the second
component $\eta _{0}$ does not necessarily depend on $u_{0}(\cdot )$. This
permits us to construct a dynamical system in this phase space and prove the
existence of global attractors. However, the transformed equation %
\eqref{eq3.3a} is a generalization of problem \eqref{eq1-1bis}, and
therefore, not every solution to equation \eqref{eq3.3a} possesses a
corresponding one to \eqref{eq1-1bis}. Notice that both problems are
equivalent if and only if the initial value $\eta _{0}$ belongs to a proper
subspace of $L_{\mu }^{2}(\mathbb{R}^{+};H_{0}^{1}(\Omega ))$, which
coincides with the domain of the distributional derivative with respecto to $%
s$, denoted by $D(\mathbf{T})$ (for more details, see \cite{G}). 
Hence, it is natural to construct a dynamical system generated by %
\eqref{eq3.3a} in the phase space $L^{2}(\Omega )\times D(\mathbf{T})$ to
prove the existence of attractors to the original problem, via the above
relationship (see \cite{C, CM, G}). Nevertheless, as far as we know, it is
not possible to prove the existence of attractors in this space unless
solutions are proved to have more regularity. Thus, in principle, we cannot
transfer the existence of attractors for system \eqref{eq3.3a} to the
original problem \eqref{eq1-1bis}.


The idea of the second alternative comes from a simple case, which was
successfully applied in \cite{CA} when the kernel is $k(t)=e^{-d_{0}t}$, $%
d_{0}>0$ (non-singular kernel). Using this method, it is proved that the
problem in \cite{CA} generates a dynamical system in the phase space $%
L_{H_{0}^{1}}^{2}$ given by the measurable functions $\varphi :(-\infty
,0]\rightarrow H_{0}^{1}(\Omega )$, such that $\int_{-\infty }^{0}e^{\gamma
s}\Vert \varphi (s)\Vert _{H_{0}^{1}}^{2}ds<+\infty $, for certain $\gamma
>0 $. Under the construction of this phase space, there exists a global
attractor to this problem (in fact, the problem in \cite{CA} is
non-autonomous and the attractor is of pullback type). Notice that, for this
kind of delay problems, in which the initial value at zero may not be
related to the values for negative times, the standard and more appropriate
phase space to construct a dynamical system is the cartesian product $%
L^{2}(\Omega )\times L_{H_{0}^{1}}^{2}$ (see \cite{CR} for more details). In
such a way, for any initial values $u_{0}\in L^{2}(\Omega )$ and $\varphi
\in L_{H_{0}^{1}}^{2}$, there exists a unique solution to the following
problem (we set $\tau =0$ since the problem is autonomous), 
\begin{equation}
\begin{cases}
c_{0}\frac{\partial u}{\partial t}-k_{0}\Delta u-\displaystyle\int_{-\infty
}^{t}k(t-s)\Delta u(s)ds+f(u)=g, \\ 
u(x,t)=0, \\ 
u(x,0)=u_{0}(x), \\ 
u(x,t)=\varphi (x,t), \\ 
\end{cases}%
\begin{aligned} &\mbox{in}~~\Omega\times (0,\infty),\\
&\mbox{on}~\partial\Omega\times\mathbb{R},\\ &\mbox{in}~\Omega,\\
&\mbox{in}~\Omega\times(-\infty,0).\\ \end{aligned}   \label{eq3.1tris}
\end{equation}%
According to the regularity of solutions to the above equation, one can
define a dynamical system $S(t):L^{2}(\Omega )\times
L_{H_{0}^{1}}^{2}\rightarrow L^{2}(\Omega )\times L_{H_{0}^{1}}^{2}$ by the
relation 
\begin{equation*}
S(t)(u_{0},\varphi ):=(u(t;0,u_{0},\varphi ),u_{t}(\cdot ;0,u_{0},\varphi )),
\end{equation*}%
where $u(\cdot ;0,u_{0},\varphi )$ denotes the solution of problem %
\eqref{eq3.1tris} (see \cite{CR} for more details on this set-up). We
emphasize that the two components of the dynamical system are the current
state of the solution and the past history up to present, respectively, what
is more sensible in a problem with delays or memory. By using this
framework, the method in \cite{CA} can be successfully applied to prove the
existence of attractors to problem \eqref{eq3.1tris} when $k$ is of
exponential type. However, this exponential behavior may be a big
restriction on the kernel $k$, consequently, on the function $\mu $, since
in many real situations the latter often has singularities, for instance $%
k(t)=e^{-d_{0}t}t^{-\alpha },\alpha \in (0,1)$. Therefore, it is interesting
to design a technique which allows us to handle the cases with this kind of
singular kernels within the context of the phase space $L^{2}(\Omega )\times
L_{H_{0}^{1}}^{2}$. We will obtain this result as a consequence of the
analysis performed in this paper even for the more general case of non-local
problems as described below.

Let us recall now that amongst many interesting results concerning non-local
differential equations, we mention the pioneering work \cite{F}, in which a
model of single-species dynamics incorporating non-local effects was
analyzed, comparing with the standard approach to model a single-species
domain $\Omega$ of ``Kolmogorov" type, 
\begin{equation*}
\partial_t u=\Delta u+\lambda u g(u),\qquad \mbox{in}\quad \Omega, ~~t>0.
\end{equation*}
Taking into account the following two natural assumptions: (i) a population
in which individuals compete for a shared rapidly equilibrate resource; (ii)
a population in which individuals communicate either visually or by chemical
means, then the most straightforward way of introducing non-local effects is
to consider, instead of $g(u)$, a ``crowding" effect of the form $g(u,\bar{u}%
)$, where 
\begin{equation*}
\bar{u}(x,t)=\int_{\Omega}G(x,y)u(y,t)dy,
\end{equation*}
and $G(x,y)$ is some reasonable kernel. Reasoning in a heuristic way, Chipot
et al. \cite{C3} studied the behavior of a population of bacteria with
non-local term $a(\int_{\Omega}u)$ in a container. Later, Chipot et al. (cf. 
\cite{C1, C2}) extended this term to a general non-local operator $a(l(u))$,
where $l\in \mathcal{L}(L^2(\Omega);\mathbb{R})$, for instance, if $g\in
L^2(\Omega)$, 
\begin{equation*}
l(u)=l_g(u)=\int_{\Omega}g(x)u(x)dx.
\end{equation*}

Motivated by these works, the dynamics of the following non-autonomous
non-local partial differential equations with delay and memory was
investigated in \cite{X1} by using the Galerkin method and energy
estimations, 
\begin{alignat}{3}  \label{eq1-2}
\begin{cases}
\frac{\partial u}{\partial t}-a(l(u))\Delta u=f(u)+h(t,u_t), \\ 
u=0, \\ 
u_{\tau}(x,t)=\varphi(x,t), \\ 
\end{cases}
\begin{aligned} &\mbox{in}~~\Omega\times (\tau,\infty),\\
&\mbox{on}~\partial\Omega\times\mathbb{R},\\ &\mbox{in}~~\Omega\times
(-\rho,0],\\ \end{aligned} & & & & 
\end{alignat}
where $\Omega\subset\mathbb{R}^N$ is a bounded open set, $\tau\in\mathbb{R}$%
, the function $a\in C(\mathbb{R};\mathbb{R}^+)$ is locally Lipschitz, $f\in
C(\mathbb{R})$, $h$ contains hereditary characteristics involving delays,
and $u_t: (-\infty, 0]\rightarrow \mathbb{R}$ is a segment of the solution
given by $u_t(x,s)=u(x,t+s)$, $s\leq0$, which essentially represents the
history of the solution up to time $t$. Moreover, $0<\rho\leq \infty$, which
implies, the authors considered both cases, bounded and unbounded delays,
for this model. However, the technique applied in \cite{X1} is the same used
in \cite{CA} and, therefore, it is valid only for non-singular memory terms
of exponential kind (e.g., $k(t)=k_1e^{-d_0t}$, $k_1\in\mathbb{R}, d_0>0$),
for more details, see \cite{CA}. Whereas, this technique fails to deal with
various important models with memory, whose kernels have singularities.%



Consequently, very recently, a new model has been considered related to long
time memory differential equations containing non-local diffusion, 
\begin{equation}
\begin{cases}
\frac{\partial u}{\partial t}-a(l(u))\Delta u-\displaystyle\int_{-\infty
}^{t}k(t-s)\Delta u(s)ds+f(u)=g, \\ 
u(x,t)=0, \\ 
u(t+\tau )=\varphi (t), \\ 
\end{cases}%
\begin{aligned} &\mbox{in}~~\Omega\times (\tau,\infty),\\
&\mbox{on}~\partial\Omega\times\mathbb{R},\\ &\mbox{in}~~\Omega\times
(-\infty,0],\\ \end{aligned}  \label{eq3.1}
\end{equation}%
where $\Omega \subset \mathbb{R}^{N}$ is a bounded domain with regular
boundary, the function $a\in C(\mathbb{R};\mathbb{R}^{+})$ satisfies 
\begin{equation}
0<m\leq a(r),\qquad \forall r\in \mathbb{R}.  \label{eq2.2A}
\end{equation}%
$k:\mathbb{R}^{+}\rightarrow \mathbb{R}$ is the memory kernel, with or
without singularities, whose properties will be specified later, $g\in
L^{2}(\Omega )$ which is independent of time. Notice that, thanks to a
change of variable, the long time memory term in problem \eqref{eq3.1} can
be interpreted as an infinite delay term, 
\begin{equation}
h(u_{t}):=\int_{-\infty }^{0}k(-s)\Delta u_{t}(x,s)ds=\int_{-\infty
}^{0}k(-s)\Delta u(x,t+s)ds=\int_{-\infty }^{t}k(t-s)\Delta u(x,s)ds.
\label{eq1}
\end{equation}%
%
%
%
%
Obviously, our model is an autonomous non-local partial differential
equation. The authors first proved in \cite{dafermos} the existence and
uniqueness of solutions to \eqref{eq3.1} by using the Dafermos
transformation. Next, they constructed an autonomous dynamical system in the
phase space $L^{2}(\Omega )\times L_{\mu }^{2}(\mathbb{R}^{+};H_{0}^{1}(%
\Omega ))$ and proved the existence of a global attractor in this space. As
in the local heat equation case mentioned above, the same lack of enough
regularity does not allow us to obtain an appropriate attractor for the
original problem \eqref{eq3.1} in the phase space $L^{2}(\Omega )\times
L_{H_{0}^{1}}^{2}$. Therefore, our objective is to overcome this difficulty
and we succeeded by proceeding in the following way: Consider problem %
\eqref{eq3.1} with initial values $u(\tau )=u_{0}$ and $u(t+\tau )=\varphi
(t)$ for $t<0$, where $(u_{0},\varphi )\in L^{2}(\Omega )\times
L_{H_{0}^{1}}^{2}$. Thus, for those kernels $\mu (\cdot )$ which guarantee
that, when $\varphi \in L_{H_{0}^{1}}^{2}$ the corresponding $\eta _{\varphi
}$, defined by $\eta _{\varphi }(s)=\int_{-s}^{0}\varphi (r)\,dr,\ (s>0)$
belongs to the space $L_{\mu }^{2}(\mathbb{R}^{+};H_{0}^{1}(\Omega ))$, we
can perform the Dafermos transformation and obtain the initial value problem
which was already analyzed in \cite{dafermos}, and consequently we have the
existence, uniqueness and regularity of solutions in a straightforward way.
Thanks to this result, we are able to construct the dynamical system in the
phase space $L^{2}(\Omega )\times L_{H_{0}^{1}}^{2}$ with the help of some
additional technical results. The existence of global attractor is then
proved by first showing the existence of a bounded absorbing set and the
proof of the asymptotic compactness property which requires an appropriate
adaptation of the technique used in \cite{CA}. These results proved in the
non-local problem \eqref{eq3.1} improve and complete the ones in \cite{CA}
by simply assuming that $a(\cdot )$ is a constant, and also improve the
previous literature on the local case (see, e.g., \cite{G, P2, C}), where it
is only provided the existence of attractors for the transformed equation %
\eqref{eq3.3a} but not for the original one \eqref{eq1-1}.

The content of this paper is as follows: In Section~\ref{s2}, we recall some
preliminaries, notations and the framework in which we will carry out our
analysis. Section~\ref{s3} is devoted to proving the main results of our
paper. First, we state the existence and uniqueness of solutions of our
problem by rewriting it as an equivalent one thanks to the Dafermos
transformation. The transformed problem has already been analyzed in \cite%
{dafermos}, whence our result follows immediately. However, as some
estimations we need for the subsequent results are based on the ones in the
proof of this existence theorem, we have included the complete proof in the
Appendix (at the end of the paper). Next, we prove that our model generates
an autonomous dynamical system in the phase space $L^2(\Omega)\times
L^2_{H^1_0}$. Eventually, the existence of a global attractor for the
dynamical system is proved by working directly on our model with memory,
instead of using any result already proved in \cite{dafermos} for the
transformed problem.

\section{Well-posedness to a non-local differential equation with memory}

\label{s2} The following non-local differential equation associated with
singular memory will be investigated, 
\begin{equation}  \label{eq0}
\begin{cases}
\frac{\partial u}{\partial t}-a(l(u))\Delta u-\displaystyle%
\int_{-\infty}^{t}k(t-s)\Delta u(x,s)ds+f(u)=g(x,t), \\ 
u(x,t)=0, \\ 
u(x,0)=u_0(x), \\ 
u(x,t+\tau)=\phi(x,t),%
\end{cases}
\begin{aligned} &\mbox{in}~~\Omega\times (\tau,\infty),\\
&\mbox{on}~\partial\Omega\times\mathbb{R},\\ &\mbox{in}~~\Omega\\
&\mbox{in}~~\Omega\times (-\infty,0],\\ \end{aligned}
\end{equation}
where $\Omega\subset\mathbb{R}^{N}$ is a fixed bounded domain with regular
boundary. The function $a\in C(\mathbb{R};\mathbb{R}^{+})$ satisfies 
\begin{equation}
0<m\leq a(r),\qquad\forall r\in\mathbb{R},  \label{eq2.2}
\end{equation}
$k:\mathbb{R}^{+}=(0,+\infty)\rightarrow\mathbb{R}$ is the memory kernel,
whose properties will be specified later. The initial values are $u_0\in
L^2(\Omega)$ and $\phi\in L^2_V$ (see Section \ref{s2.2} below).

Let us define the new variables 
\begin{equation*}
u^{t}(x,s)=u(x,t-s),\qquad s\geq0,
\end{equation*}
and 
\begin{equation}
\eta^{t}(x,s)=\int_{0}^{s}u^{t}(x,r)dr=\int_{t-s}^{t}u(x,r)dr,\qquad s\geq0.
\label{eq3.2}
\end{equation}
Assuming $k(\infty)=0$, a change of variable and a formal integration by
parts yield 
\begin{equation*}
\int_{-\infty}^{t}k(t-s)\Delta
u(s)ds=-\int_{0}^{\infty}k^{\prime}(s)\Delta\eta^{t}(s)ds,
\end{equation*}
here and in the sequel, the prime denotes derivation with respect to
variable $s$. Setting 
\begin{equation}
\mu(s)=-k^{\prime}(s),  \label{eq3.01}
\end{equation}
the above choice of variable leads to the following non-delay system, 
\begin{equation}
\begin{cases}
\frac{\partial u}{\partial t}-a(l(u))\Delta u-\displaystyle%
\int_{0}^{\infty}\mu(s)\Delta\eta^{t}(s)ds+f(u)=g(x,t), \\ 
\frac{\partial}{\partial t}\eta^{t}(s)=u-\frac{\partial}{\partial s}%
\eta^{t}(s), \\ 
u(x,t)=\eta^{t}(x,s)=0, \\ 
u(x,\tau)=u_{0}(x), \\ 
\eta^{\tau}(x,s)=\eta_{0}(x,s),%
\end{cases}
\begin{aligned} &\mbox{in}~~\Omega\times (\tau,\infty),\\
&\mbox{in}~~\Omega\times (\tau,\infty)\times\mathbb{R}^+,\\
&\mbox{on}~~\partial \Omega\times \mathbb{R}\times\mathbb{R}^+,\\
&\mbox{in}~~ \Omega,\\ &\mbox{in}~~\Omega\times\mathbb{R}^+, \end{aligned}
\label{eq3.3}
\end{equation}
where, by the definition of $\eta^{t}(x,s)$ (see (\ref{eq3.2})), it
obviously follows 
\begin{equation}
\eta^{\tau}(x,s)=\int_{\tau-s}^{\tau}u(x,r)dr=\int_{-s}^{0}\phi(x,r)dr:=\eta
_{0}(x,s),  \label{eq2.5}
\end{equation}
which is the initial integrated past history of $u$ with vanishing boundary.

It is worth emphasizing that we will consider solutions of our problems in
the weak (variational) sense.

\indent{\subsection{Assumptions}}

In our analysis, we shall suppose the nonlinear term $f:\mathbb{R}\rightarrow%
\mathbb{R}$ is a polynomial of odd degree with positive leading coefficient, 
\begin{equation}
f(u)=\sum_{k=1}^{2p}f_{2p-k}u^{k-1},\quad ~~p\in\mathbb{N}.  \label{eq3.4}
\end{equation}
This situation can be extended, without any additional difficulties, to a
more general function satisfying suitable assumptions (see, for instance, 
\cite{C}).

In view of the evolution problem (\ref{eq3.3}), the variable $\mu$ is
required to verify the following hypotheses:

\begin{enumerate}
\item[$(h_{1})$] $\mu\in C^{1}(\mathbb{R}^{+})\cap L^{1}(\mathbb{R}%
^{+}),\quad\mu(s)\geq0,\quad\mu^{\prime}(s)\leq0$,\quad$\forall s\in\mathbb{R%
}^{+}$;

\item[$(h_{2})$] $\mu^{\prime}(s)+\delta\mu(s)\leq0$,\quad$\forall s\in%
\mathbb{R}^{+}$, ~~for some~ $\delta>0$.
\end{enumerate}


\begin{remark}
\label{rem3.2}

\begin{enumerate}
\item It is straightforward to check that conditions $(h_1)$-$(h_2)$ are
fulfilled by singular kernels given by 
\begin{equation*}
\mu(t)=e^{-\delta t}t^{-\alpha},\ t>0,
\end{equation*}
for $\delta>0$ and $\alpha\in(0,1)$.

\item Restriction $(h_{1})$ is equivalent to assuming $k(\cdot )$ is a
non-negative, non-increasing, bounded, convex function of class $C^{2}$
converging at infinity. Moreover, from $(h_{1})$ and $k(t)\underset{%
t\rightarrow +\infty }{\rightarrow }0$ it easily follows that 
\begin{equation*}
k(0)=\int_{0}^{\infty }\mu (s)ds
\end{equation*}
is finite and non-negative.

\item Assumption $(h_{2})$ implies that $\mu (s)$ decays exponentially.
Indeed, a simple integration shows that%
\begin{equation*}
\mu (s)\leq \mu (t)e^{-\delta (s-t)}\text{ for }0<t<s.
\end{equation*}

Also, this condition allows the memory kernel $k(\cdot )$ to have a
singularity at $t=0$, which coincides with the intention to study problem (%
\ref{eq3.3}). 

\item Since it is assumed that $\lim_{t\rightarrow +\infty }k(t)=0$, we have%
\begin{equation*}
k(t)=-\int_{t}^{\infty }k^{\prime }\left( s\right) ds=\int_{t}^{\infty }\mu
\left( s\right) ds\leq M_{1},\ \ \forall t\geq 0,
\end{equation*}%
for some positive constant $M_{1}.$ Also, 
\begin{equation*}
k(t)=\int_{t}^{\infty }\mu \left( s\right) ds\leq \mu (t)\int_{t}^{\infty
}e^{-\delta (s-t)}ds=\frac{\mu (t)}{\delta }.
\end{equation*}
\end{enumerate}
\end{remark}

\subsection{Notations}

\label{s2.2} {Let $\Omega$ be a fixed bounded domain in $\mathbb{R}^N$. On
this set, we introduce the Lebesgue space $L^p(\Omega)$, where $1\leq p\leq
\infty$. Besides, $W^{1,p}(\Omega)$ is the subspace of $L^p( \Omega)$
consisting of functions such that the first order weak derivative belongs to 
$L^p(\Omega)$. In this paper, $L^2(\Omega)$ is denoted by $H$, $%
H_0^1(\Omega) $ is denoted by $V$ and $H^{-1}(\Omega)$ is denoted by $V^*$.
The norms in $H $, $V$ and $V^*$ will be denoted by $|\cdot|$, $\|\cdot\|$
and $\| \cdot\|_*$, respectively. }

In view of system (\ref{eq3.3}) and $(h_1)$, we need to introduce some
additional notations before proving our main theorems. Let $L^2_{\mu}(%
\mathbb{R}^+;H)$ be a Hilbert space of functions $w:\mathbb{R}^+\rightarrow
H $ endowed with the inner product, 
\begin{equation*}
(w_1,w_2)_{\mu}=\int_{0}^{\infty}\mu(s)(w_1(s),w_2(s))ds,
\end{equation*}
and let $|\cdot|_{\mu}$ denote the corresponding norm. In a similar way, we
introduce the inner products $((\cdot,\cdot))_{\mu}$, $(((\cdot,\cdot)))_{%
\mu}$ and relative norms $\|\cdot\|_{\mu}$, $|||\cdot|||_{\mu}$ on $%
L^2_{\mu}(\mathbb{R}^+;V)$, $L^2_{\mu}(\mathbb{R}^+;V\cap H^2(\Omega))$
respectively. It follows then that 
\begin{equation*}
((\cdot,\cdot))_{\mu}=(\nabla\cdot,\nabla\cdot)_{\mu},~~\mbox{and}%
~~(((\cdot,\cdot)))_{\mu}=(\Delta\cdot,\Delta\cdot)_{\mu}.
\end{equation*}
We also define the Hilbert spaces 
\begin{equation*}
\mathcal{H}=H\times L^2_{\mu}(\mathbb{R}^+;V),
\end{equation*}
and 
\begin{equation*}
\mathcal{V}=V\times L^2_{\mu}(\mathbb{R}^+;V\cap H^2(\Omega)),
\end{equation*}
which are respectively endowed with inner products 
\begin{equation*}
(w_1,w_2)_{\mathcal{H}}=(w_1,w_2)+((w_1,w_2))_{\mu},
\end{equation*}
and 
\begin{equation*}
(w_1,w_2)_{\mathcal{V}}=((w_1,w_2))+(((w_1,w_2)))_{\mu},
\end{equation*}
where $w_i \in \mathcal{H}$ or $\mathcal{V}$ $(i=1,2)$ and usual norms.

At last, with standard notations, $\mathcal{D}(I;X)$ is the space of
infinitely differentiable $X$-valued functions with compact support in $%
I\subset \mathbb{R}$, whose dual space is the distribution space on $I$ with
values in $X^*$ (dual of $X$), denoted by $\mathcal{D}^{\prime}(I;X^{*})$.
For convenience, we define $L_{V}^{2}$ the space of functions $u\left(\text{%
\textperiodcentered}\right)$ satisfying%
\begin{equation*}
\int_{-\infty}^{0}e^{\gamma s}\left\Vert u\left( s\right) \right\Vert
^{2}ds<\infty,
\end{equation*}
where $0<\gamma<\min\{m\lambda_{1},\delta\}$ and $\delta$ comes from $(h_2)$%
. 

\section{Main results}

\label{s3} Let us start by proving a technical result which will be crucial
to our analysis. To this end, we define the linear operator $\mathcal{J}:
L^2_V\rightarrow L^2_\mu(\mathbb{R}^+;V)$ by 
\begin{equation}  \label{j}
(\mathcal{J}\phi)(s)=\int_{-s}^0\phi(r)\,dr,\quad s\in\mathbb{R}^+.
\end{equation}
Then we have the following result.

\begin{lemma}
\label{lemma3.1} Assume $(h_1)$-$(h_2)$ hold. Then, the operator $\mathcal{J}
$ defined by \eqref{j} is a linear and continuous mapping. In particular,
there exists a positive constant $K_\mu$ such that, for any $\phi\in L_V^2$,
it holds 
\begin{equation}  \label{Kmu}
\|\mathcal{J}\phi\|^2_{L_{\mu}^2(\mathbb{R}^+;V)}\leq
K_\mu\|\phi\|^2_{L^2_V}.
\end{equation}
\end{lemma}

\noindent \textbf{Proof.} The linearity of $\mathcal{J}$ is obvious, we only
need to prove it is well defined and bounded. Indeed, taking into account
the fact that $\phi\in L^2_V$, $(h_1)$-$(h_2)$ and \eqref{j}, we have 
\begin{equation*}
\begin{split}
\|\mathcal{J}\phi\|^2_{L_{\mu}^2(\mathbb{R}^+;V)}&=\int_0^{\infty}\mu(s)%
\left\|\int_{-s}^0\phi(r)dr\right\|^2ds \\
&=\int_0^{1}\mu(s)\left\|\int_{-s}^0\phi(r)dr\right\|^2ds+\int_{1}^{\infty}%
\mu(s)\left\|\int_{-s}^0\phi(r)dr\right\|^2ds \\
&\leq \int_0^{1}s\mu(s)\int_{-s}^0\|\phi(r)\|^2drds+
\mu(1)\int_1^{\infty}e^{-\delta (s-1)}\left\|\int_{-s}^0\phi(r)dr\right\|^2ds
\\[1.0ex]
&\leq
\int_{-1}^0\|\phi(r)\|^2\int_{-r}^{1}s\mu(s)dsdr+\mu(1)e^{\delta}\int_0^{%
\infty}e^{-\delta s}s\int_{-s}^0\|\phi(r)\|^2drds \\
&\leq
\int_0^{1}s\mu(s)ds\int_{-1}^0\|\phi(r)\|^2dr+\mu(1)e^{\delta}\int_{-%
\infty}^0e^{\gamma r}\|\phi(r)\|^2\int_{-r}^{\infty}se^{-\gamma r}e^{-\delta
s}dsdr \\[1.0ex]
&\leq \int_0^{1}\mu(s)ds\int_{-1}^0e^{-\gamma r}e^{\gamma r}\|\phi(r)\|^2dr
\\
&\quad+ \mu(1)e^{\delta}\int_{-\infty}^0 e^{\gamma
r}\|\phi(r)\|^2\int_{-r}^{\infty} se^{\gamma s} e^{-\delta s}dsdr \\
&\leq \left(e^{\gamma }\int_0^{1} \mu(s)ds+
\mu(1)e^{\delta}(\gamma-\delta)^{-2} \right)\|\phi\|^2_{L_V^2}.
\end{split}%
\end{equation*}
Denoting $K_{\mu}=e^{\gamma}\int_0^{1} \mu(s)ds+
\mu(1)e^{\delta}(\gamma-\delta)^{-2} $, the proof is finished. $\Box$

\begin{remark}
Notice that when we fix an initial value $\phi\in L_V^2$ for problem %
\eqref{eq0}, then the corresponding initial value for the second component
of problem \eqref{eq3.3} becomes $\eta_0:=\mathcal{J}\phi$, which belongs to 
$L^2_{\mu}(\mathbb{R}^+;V)$ thanks to Lemma \ref{lemma3.1}. 
\end{remark}

Before stating the existence, uniqueness and regularity of solution to our
problem \eqref{eq0}, we first recall a general result proved in \cite%
{dafermos} for problem \eqref{eq3.3} with general initial data in $H\times
L_{\mu }^{2}(\mathbb{R}^{+};V)$. Let us denote 
\begin{equation*}
z(t)=(u(t),\eta ^{t})\qquad \mbox{and}\qquad z_{0}=(u_{0},\eta _{0}).
\end{equation*}%
Set 
\begin{equation*}
\mathcal{L}z=\left( a(l(u))\Delta u+\int_{0}^{\infty }\mu (s)\Delta \eta
(s)ds,~~u-\eta _{s}\right) ,
\end{equation*}%
and 
\begin{equation*}
\mathcal{G}(z)=(-f(u)+g,~~0).
\end{equation*}%
Then problem (\ref{eq3.3}) can be written in the following compact form, 
\begin{equation}
\begin{cases}
z_{t}=\mathcal{L}z+\mathcal{G}(z), \\ 
z(x,t)=0, \\ 
z(x,\tau )=z_{0},%
\end{cases}%
\begin{aligned} &\qquad \mbox{in}~~ \Omega\times (\tau,\infty),\\ &\qquad
\mbox{on}~~\partial\Omega\times(\tau,\infty),\\ &\qquad \mbox{in}~~\Omega.
\end{aligned}  \label{eq3.5}
\end{equation}

Now we have the following result.

\begin{theorem}[\protect\cite{dafermos}]
\label{thm3.1} Suppose \eqref{eq2.2}, \eqref{eq3.4} and $(h_1)$-$(h_2)$ hold
true, also let $g\in H$. In addition, assume that $a(\cdot)$ is locally
Lipschitz, and there exists a positive constant $\tilde{m}$ such that, 
\begin{equation}  \label{m}
a(s)\leq \tilde{m},\qquad \forall s\in\mathbb{R}.
\end{equation}
Then:

\begin{enumerate}
\item[$(i)$] For any $z_{0}\in \mathcal{H}$, there exists a unique solution $%
z(\cdot )=(u(\cdot ),\eta ^{\cdot })$ to problem \eqref{eq3.5} which
satisfies 
\begin{equation*}
\begin{split}
& u(\cdot )\in L^{\infty }(\tau ,T;H)\cap L^{2}(\tau ,T;V)\cap L^{2p}(\tau
,T;L^{2p}(\Omega )),\qquad \forall T>\tau , \\
& \eta ^{\cdot }\in L^{\infty }(\tau ,T;L_{\mu }^{2}(\mathbb{R}%
^{+};V)),\qquad \forall T>\tau .
\end{split}%
\end{equation*}%
Furthermore, $z(\cdot )\in C([\tau ,T];\mathcal{H})$ for every $T>\tau $,
and the mapping $F:z_{0}\in \mathcal{H}\rightarrow z(t)\in \mathcal{H}$ is
continuous for every $t\in \lbrack \tau ,T]$.

\item[$(ii)$] For any $z_{0}\in \mathcal{V}$, the unique solution $z(\cdot
)=(u(\cdot ),\eta ^{\cdot })$ to problem \eqref{eq3.5} satisfies 
\begin{equation*}
\begin{split}
& u(\cdot )\in L^{\infty }(\tau ,T;V)\cap L^{2}(\tau ,T;V\cap H^{2}(\Omega
)),\qquad \forall T>\tau , \\
& \eta ^{\cdot }\in L^{\infty }(\tau ,T;L_{\mu }^{2}(\mathbb{R}^{+};V\cap
H^{2}(\Omega ))),\qquad \forall T>\tau .
\end{split}%
\end{equation*}%
In addtion, $z(\cdot )\in C([\tau ,T];\mathcal{V})$ for every $T>\tau $.
\end{enumerate}
\end{theorem}

Based on the previous theorem, we can state now the corresponding result for
our original problem \eqref{eq0}.

\begin{theorem}
\label{Existence} Assume (\ref{eq2.2}), (\ref{eq3.4}), and $(h_{1})$-$(h_{2})
$ hold. Let $a(\cdot )$ be locally Lipschitz satisfying \eqref{m}, 
\begin{equation*}
g\in H,\quad u_{0}\in H\quad \mbox{and}\quad \phi \in L_{V}^{2}.
\end{equation*}%
Then, there exists a unique function $z(\cdot )=(u(\cdot ),\eta ^{\cdot })$
satisfying 
\begin{equation*}
\begin{split}
& u(\cdot )\in L^{\infty }(\tau ,T;H)\cap L^{2}(\tau ,T;V)\cap L^{2p}(\tau
,T;L^{2p}(\Omega )),\qquad \forall T>\tau , \\
& \eta ^{\cdot }\in L^{\infty }(\tau ,T;L_{\mu }^{2}(\mathbb{R}%
^{+};V)),\qquad \forall T>\tau ,
\end{split}%
\end{equation*}%
such that 
\begin{equation*}
\partial _{t}z=\mathcal{L}z+\mathcal{G}(z)
\end{equation*}%
in the weak sense, and 
\begin{equation*}
z|_{t=\tau }=(u_{0},\mathcal{J}\phi ).
\end{equation*}%
Furthermore, for every $t\in \lbrack \tau ,T]$, $z_{0}\mapsto $ $z(t)~$is a
continuous mapping from $\mathcal{H}$\ into $\mathcal{H}$. 

If we also assume that $u_{0}\in V$, $\phi \in L_{V\cap H^{2}(\Omega )}^{2}$%
, then 
\begin{equation*}
\begin{split}
& u\in L^{\infty }(\tau ,T;V)\cap L^{2}(\tau ,T;V\cap H^{2}(\Omega )),\qquad
\forall T>\tau , \\
& \eta ^{\cdot }\in L^{\infty }(\tau ,T;L_{\mu }^{2}(\mathbb{R}^{+};V\cap
H^{2}(\Omega ))),\qquad \forall T>\tau ,
\end{split}%
\end{equation*}%
and for each $t\in \lbrack \tau ,T]$, $z_{0}\mapsto $ $z(t)~$is a continuous
mapping from $\mathcal{V}$ into $\mathcal{V}.$.
\end{theorem}

\indent {\bf Proof.} Thanks to Lemma~\ref{lemma3.1}, we obtain $\mathcal{J}%
\phi\in L_{\mu}^{2}(\mathbb{R}^{+};V)$ since $\phi\in L^2_V$. Therefore, the
first statement of Theorem~\ref{Existence} holds by applying (i) in Theorem~%
\ref{thm3.1} with initial value $z_0=(u_0,\mathcal{J}\phi)$. If, in
addition, we assume that initial values $u_0\in V$ and $\phi\in L^2_{V\cap
H^2(\Omega)}$, then it is straightforward to prove that $z_0=(u_0,\mathcal{J}%
\phi)\in\mathcal{V}$ and the regularity result follows from statement (ii)
in Theorem~\ref{thm3.1}. $\Box$

\begin{remark}
Although the proof of Theorem~\ref{Existence} follows directly from Theorem~%
\ref{thm3.1}, some computations, that we need in the sequel, are based on
some estimations carried out in the proof. For this reason, we have included
the complete proof of Theorem~\ref{Existence} as an Appendix, so that the
paper is self-contained and easier to read.
\end{remark}

The following lemma allows us to show rigorously that in fact the function $%
u $ satisfies the equation in (\ref{eq0}).

\begin{lemma}
\label{lemEquivInt}Let conditions (h1)-(h2) hold. If $u\in L_{V}^{2}$, then $%
\eta(s)=\int_{-s}^{0}u\left( r\right) dr$ belongs to $L_{\mu}^{2}(\mathbb{R}%
^{+};V)$ and%
\begin{equation}
\int_{0}^{\infty}\mu(s)\Delta\eta(s)ds=\int_{-\infty}^{0}k(-s)\Delta u(s)ds.
\label{EquivIntegrals}
\end{equation}
\end{lemma}

\textbf{Proof. }The fact that $\eta \in L_{\mu }^{2}(\mathbb{R}^{+};V)$ is
given by Lemma \ref{lemma3.1}. From the arguments in \cite[pp-174-175]%
{Gajewski}, it follows the existence of a sequence of functions $u_{n}\left( 
\text{\textperiodcentered }\right) \in C^{1}((-\infty ,0],V)\cap L_{V}^{2}$
such that 
\begin{equation*}
u_{n}\rightarrow u\text{ in }L_{V}^{2}.
\end{equation*}%
First, we will show that $u_{n},\eta _{n}$, where $\eta
_{n}(s)=\int_{-s}^{0}u_{n}\left( r\right) dr$, satisfy (\ref{EquivIntegrals}%
). For any $w\in V$, we have%
\begin{align*}
\left\langle \int_{0}^{\infty }\mu (s)\Delta \eta _{n}(s)ds,w\right\rangle &
=\int_{0}^{\infty }\mu (s)\left\langle \Delta \eta _{n}(s),w\right\rangle
ds=\int_{0}^{\infty }k^{\prime }\left( s\right) \left( \nabla \eta
_{n}(s),\nabla w\right) ds \\
& =\int_{0}^{\infty }k^{\prime }\left( s\right) \left( \nabla
\int_{-s}^{0}u_{n}(r)dr,\nabla w\right) ds \\
& =\int_{0}^{\infty }k^{\prime }\left( s\right) \int_{-s}^{0}\left( \nabla
u_{n}(r),\nabla w\right) drds \\
& =-\int_{0}^{\infty }k(s)\left( \nabla u_{n}(-s),\nabla w\right)
ds+\lim_{s\rightarrow \infty }\ k(s)\int_{-s}^{0}\left( \nabla
u_{n}(r),\nabla w\right) dr \\
& ~~-\lim_{s\rightarrow 0}\ k(s)\int_{-s}^{0}\left( \nabla u_{n}(r),\nabla
w\right) dr.
\end{align*}%
Let us check that the last two limits of the above equality are equal to $0$%
. By Remark \ref{rem3.2}, we derive%
\begin{equation*}
k(s)e^{\gamma s}\leq \frac{\mu (1)}{\delta }e^{\delta }e^{(\gamma -\delta
)s},~~~\text{ for any }s\geq 1.
\end{equation*}%
Hence, $\gamma <\delta $ implies%
\begin{align*}
\left\vert k(s)\int_{-s}^{0}\left( \nabla u_{n}(r),\nabla w\right)
dr\right\vert & \leq k(s)e^{\gamma s}\left\Vert w\right\Vert
\int_{-s}^{0}e^{\gamma r}\left\Vert u_{n}(r)\right\Vert dr \\
& \leq \frac{k(s)e^{\gamma s}\left\Vert w\right\Vert }{2}\left(
\int_{-\infty }^{0}e^{\gamma r}\left\Vert u_{n}(r)\right\Vert ^{2}dr+\frac{1%
}{\gamma }\right) \leq C_{1}e^{(\gamma -\delta )s}\underset{s\rightarrow
\infty }{\rightarrow }0.
\end{align*}%
Also, from $k\left( s\right) \underset{s\rightarrow 0}{\rightarrow }%
\int_{0}^{\infty }\mu \left( r\right) dr$ and $u_{n}\in L_{V}^{2}$, it
follows that the second limit is $0$ as well. Hence,%
\begin{equation*}
\left\langle \int_{0}^{\infty }\mu (s)\Delta \eta _{n}(s)ds,w\right\rangle
=-\int_{0}^{\infty }k(s)\left( \nabla u_{n}(-s),\nabla w\right)
ds=\left\langle \int_{-\infty }^{0}k(-s)\Delta u_{n}(s)ds,w\right\rangle .
\end{equation*}%
This proves (\ref{EquivIntegrals}) for $u_{n}$.

Furthermore, for any $w\in V$, we infer%
\begin{align*}
& \left\vert \left\langle \int_{-\infty }^{0}k(-s)(\Delta u_{n}(s)-\Delta
u(s))ds,w\right\rangle \right\vert  \\
& =\left\vert \int_{-\infty }^{0}k(-s)\left\langle \Delta u_{n}(s)-\Delta
u(s),w\right\rangle ds\right\vert  \\
& \leq \left\Vert w\right\Vert \left( C_{2}\int_{-1}^{0}\left\Vert
u_{n}(s)-u(s)\right\Vert ds+\frac{\mu (1)}{\delta }e^{\delta }\int_{-\infty
}^{-1}e^{\delta s}\left\Vert u_{n}(s)-u(s)\right\Vert ds\right)  \\
& \leq C_{3}\left( \left( \int_{-1}^{0}\left\Vert u_{n}(s)-u(s)\right\Vert
^{2}ds\right) ^{\frac{1}{2}}+\left( \int_{-\infty }^{-1}e^{\gamma
s}\left\Vert u_{n}(s)-u(s)\right\Vert ^{2}ds\right) ^{\frac{1}{2}}\right) 
\underset{n\rightarrow \infty }{\rightarrow }0,
\end{align*}%
and Lemma \ref{lemma3.1} implies%
\begin{align*}
& \left\vert \left\langle \int_{0}^{\infty }\mu (s)(\Delta \eta
_{n}(s)-\Delta \eta (s))ds,w\right\rangle \right\vert  \\
& =\left\vert \int_{0}^{\infty }\mu (s)\left\langle \Delta \eta
_{n}(s)-\Delta \eta (s),w\right\rangle ds\right\vert \leq \left\Vert
w\right\Vert \int_{0}^{\infty }\mu (s)\left\Vert \eta _{n}(s)-\eta
(s)\right\Vert ds \\
& \leq \left\Vert w\right\Vert \left( \int_{0}^{\infty }\mu (s)ds\right) ^{%
\frac{1}{2}}\left( \int_{0}^{\infty }\mu (s)\left\Vert \eta _{n}(s)-\eta
(s)\right\Vert ^{2}ds\right) ^{\frac{1}{2}}\leq C_{4}\left\Vert
u_{n}-u\right\Vert _{L_{V}^{2}}\underset{n\rightarrow \infty }{\rightarrow }%
0.
\end{align*}%
By these convergences we deduce (\ref{EquivIntegrals}). The proof of this
lemma is complete. $\Box $

\bigskip

Lemma \ref{lemEquivInt} implies that the solution given in Theorem \ref%
{Existence} is in fact the unique weak solution to problem (\ref{eq0}).

\begin{corollary}
\label{CorExistSol}Assume the conditions of Theorem \ref{Existence}. Then if
the function $(u,\eta )$ {is} the unique weak solution to problem (\ref%
{eq3.5}) corresponding to the initial values $u_{0}\in H$ and $\varphi \in
L_{V}^{2}$, then $u\left( \text{\textperiodcentered }\right) $ is the unique
weak solution to problem (\ref{eq0}).
\end{corollary}

\bigskip

In what follows, we construct the dynamical system generated by \eqref{eq0}
assuming that $g$ does not depend on $t$, which makes our problem be
autonomous. Thus, the theory of autonomous dynamical systems is appropriate
to carry out the analysis of the global asymptotic behavior. We emphasize
that the non-autonomous case can also be studied by exploiting the theory of
non-autonomous dynamical systems (either the theory of pullback attractors
or the uniform attractors one). The autonomous framework is concerned with
the phase space 
\begin{equation*}
X=H\times L_{V}^{2},
\end{equation*}%
endowed with the norm 
\begin{equation*}
\Vert (w_{1},w_{2})\Vert _{X}^{2}=|w_{1}|^{2}+\Vert w_{2}\Vert
_{L_{V}^{2}}^{2}.
\end{equation*}%
Then, thanks to Theorem \ref{Existence}, we can define a semigroup $S:%
\mathbb{R}^{+}\times X\rightarrow X$ by 
\begin{equation*}
S(t)\left( u_{0},\phi \right) =(u(t;0,(u_{0},\mathcal{J}\phi )),u_{t}(\cdot
;0,(u_{0},\mathcal{J}\phi ))),
\end{equation*}%
where $\left( u(\cdot ;0,(u_{0},\mathcal{J}\phi )),\eta ^{\cdot }\right) $
is the unique solution to problem (\ref{eq3.3}) with $u\left( 0\right)
=u_{0} $, $\eta _{0}=\mathcal{J}\phi $.

Let us first prove that the dynamical system $S$ is well defined. In what
follows, we will take $\tau=0$ since we are working on autonomous dynamical
system.

\begin{lemma}
Under assumptions of Theorem \ref{Existence}, if $\left(u_{0},\phi\right)\in
X$, then $S(t)\left(u_{0},\phi\right)\in X.$
\end{lemma}

\textbf{Proof.} Let $\left(u_{0},\phi\right)\in X$ and, for simplicity,
denote by $(u(\cdot),\eta^{\cdot})$ the solution to problem (\ref{eq3.3})
corresponding to the initial value $(u_0,\mathcal{J}\phi)$. It follows from
Theorem \ref{Existence} that the first component $u(t)$ belongs to $H$, thus
it only remains to show that the segment of solution $u_t(\cdot)$ belongs to 
$L^2_V$. Indeed, 
\begin{eqnarray*}
\int_{-\infty}^0 e^{\gamma s}\|u_t(s)\|^2\,ds&=&\int_{-\infty}^0 e^{\gamma
s}\|u(t+s)\|^2\,ds \\
&=&\int_{-\infty}^t e^{\gamma (\sigma-t)}\|u(\sigma)\|^2\,d\sigma \\
&=&e^{-\gamma t}\int_{-\infty}^t e^{\gamma\sigma}\|u(\sigma)\|^2\,d\sigma \\
&=&e^{-\gamma t}\int_{-\infty}^0
e^{\gamma\sigma}\|\phi(\sigma)\|^2\,d\sigma+ \int_{0}^t{e^{\gamma(\sigma-t)}}
\|u(\sigma)\|^2\,d\sigma \\[1.0ex]
&<&+\infty,
\end{eqnarray*}
where the above estimation holds true since $\phi\in L^2_V$ and $u\in
L^2(0,T;V)$ for all $T>0$. The proof of this lemma is complete. $\Box$

\begin{lemma}
\label{Est1} Under assumptions of Theorem~\ref{Existence}, there exist two
positive constants $K_1$ and $K_2$, such that 
\begin{equation}
\left\Vert S(t)(u_0,\phi)\right\Vert _{X}^{2}\leq K_{1}\left\Vert (u_0,\phi)
\right\Vert _{X}^{2}e^{-\gamma t}+K_{2},\qquad \forall t\geq 0,\
(u_0,\phi)\in X.  \label{IneqX}
\end{equation}
\end{lemma}

\noindent\textbf{Proof.} Let $(u_0,\phi)\in X$ and denote by $%
z(\cdot)=(u(\cdot),\eta^{\cdot})$ the solution to (\ref{eq3.3})
corresponding to the initial value $(u_0,\mathcal{J}\phi)$. Now, we multiply
the first equation in (\ref{eq3.3}) by $u\left( t\right) $ in $H$ and the
second equation in \eqref{eq3.3} by $\eta^{t}$ in $L^2_{\mu}(\mathbb{R}^+;V)$%
. Then, by the same energy estimations as in the proof of Theorem \ref%
{Existence} (see Appendix \eqref{eq611-3}), we obtain%
\begin{align*}
& \frac{d}{dt}\left\Vert z\right\Vert _{\mathcal{H}}^{2}+m\lambda
_{1}\left\vert u\right\vert ^{2}+m\left\Vert u\right\Vert
^{2}+f_{0}\left\vert u\right\vert _{2p}^{2p}+2(((\eta^{t})^{\prime},\eta
^{t}))_{\mu} \\
& \leq2a_{0}\left\vert \Omega\right\vert +\frac{2}{\sqrt{\lambda_{1}}}%
\left\vert g\right\vert \left\Vert u\right\Vert \\
& \leq2a_{0}\left\vert \Omega\right\vert +\frac{2}{m\lambda_{1}}\left\vert
g\right\vert ^{2}+\frac{m}{2}\left\Vert u\right\Vert^2.
\end{align*}
Since%
\begin{equation}
2(((\eta^{t})^{\prime},\eta^{t}))_{\mu}=-\int_{0}^{\infty}\mu^{\prime
}(s)|\nabla\eta^{t}(s)|^{2}ds\geq\delta\int_{0}^{\infty}\mu(s)|\nabla\eta
^{t}(s)|^{2}ds,  \label{IneqEta}
\end{equation}
it follows that%
\begin{equation}
\frac{d}{dt}\left\Vert z\right\Vert _{\mathcal{H}}^{2}+\gamma\left\Vert
z\right\Vert _{\mathcal{H}}^{2}+\frac{m}{2}\left\Vert u\right\Vert^2
+f_{0}\left\vert u\right\vert _{2p}^{2p}\leq K_0,  \label{Ineq0}
\end{equation}
where $K_0=2a_0|\Omega|+\frac{2}{m\lambda_1} |g|^2$ and we recall that $%
\gamma<\min\{m\lambda_1,\delta\}$. Notice that inequality \eqref{IneqEta}
has been deduced formally but can be fully justified by using mollifiers
(see \cite[p. 348]{C}). Now multiplying the above inequality by $e^{\gamma
t} $ and integrating over $\left( 0,t\right) $, neglecting the last term of
the left hand side of \eqref{Ineq0}, we obtain%
\begin{align}
&\quad \left\Vert z\left( t\right) \right\Vert _{\mathcal{H}}^{2}+\frac{m}{2}%
\int_{0}^{t}e^{-\gamma(t-s)}\left\Vert u\left( s\right) \right\Vert ^{2}ds 
\notag \\
&\leq \|z(t)\|_{\mathcal{H}}^2+\frac{m}{2}\int_{-t}^0e^{\gamma
s}\|u_t(s)\|^2ds  \notag \\
&\leq \left\Vert z_0 \right\Vert _{\mathcal{H}}^{2}e^{-\gamma t}+\frac{K_0}{%
\gamma}.  \label{Ineq1}
\end{align}
Then%
\begin{equation*}
\begin{split}
\frac{m}{2}\|u_t\|^2_{L_V^2}&=\frac{m}{2}\int_{-\infty}^0 e^{-\gamma
(t-s)}\|\phi(s)\|^2ds+\frac{m}{2}\int_{0}^te^{-\gamma(t-s)}\|u(s)\|^2ds \\
&\leq \frac{m}{2}e^{-\gamma t}\|\phi\|^2_{L_V^2}+\|(u_0,\mathcal{J}%
\phi)\|^2_{\mathcal{H}} e^{-\gamma t}+\frac{K_0}{\gamma}.
\end{split}%
\end{equation*}
In view of Lemma~\ref{lemma3.1}, we have that 
\begin{equation}
\left\Vert z_0 \right\Vert _{\mathcal{H}}^{2}\leq\left\vert u_0 \right\vert
^{2}+\|\mathcal{J}\phi\|^2_{L^2_{\mu}(\mathbb{R}^+;V)}\leq |u_0|^2+ {K_\mu}%
\|\phi\|^2_{L_V^2}.  \label{Ineq3}
\end{equation}
Hence, (\ref{Ineq1})-(\ref{Ineq3}) imply the existence of positive constants 
$K_{1}$ and $K_{2}$, such that%
\begin{equation*}
\|S(t)(u_0,\phi)\|^2_{X}:=|u(t)|^2+\|u_t\|_{L_V^2}^2\leq K_1
\left(|u_0|^2+\|\phi\|^2_{L_V^2}\right)e^{-\gamma t}+K_2.
\end{equation*}
The proof of this lemma is complete. $\Box$

From Lemma~\ref{Est1}, we immediately have the following result.

\begin{corollary}
\label{Absorbing}The ball $B_{0}=\{v\in X:\left\Vert v\right\Vert
_{X}^{2}\leq2K_{2}\}$ is absorbing for the semigroup $S$.
\end{corollary}

Now we shall prove the asymptotic compactness of the semigroup $S$. To this
end, we first state the next result.

\begin{lemma}
\label{ConvergSol}Assume the hypotheses in Theorem \ref{Existence}. Let $%
\{\left( u_{0}^{n},\phi^{n}\right) \}$ be a sequence, such that $\left(
u_{0}^{n},\phi^{n}\right) \rightarrow\left( u_{0},\phi\right) $ weakly in $X$
as $n\rightarrow \infty$. Then, $S(t)\left( u_{0}^{n},\phi^{n}\right)
=(u^n(t), u^n_t) $ fulfills:%
\begin{equation}
u^{n}\left( \cdot\right) \rightarrow u\left( \cdot\right) ~~ \text{ in }%
~~C([r,T],H) ~~\text{ for all }~~ 0<r<T;  \label{Conv1}
\end{equation}%
\begin{equation}
u^{n}(\cdot)\rightarrow u(\cdot)~~\text{ weakly in }~~L^{2}\left(0,T;V
\right) ~~\text{ for all }~~ T>0;  \label{Conv2}
\end{equation}%
\begin{equation}
u^{n}\rightarrow u~~\text{ in }~~ L^{2}\left(0, T;H \right) ~~ \text{ for
all }~~T>0;  \label{Conv2B}
\end{equation}%
\begin{equation}
\limsup_{n\rightarrow\infty}\ \left\Vert u_t^{n}-u_t\right\Vert
_{L_V^2}^{2}\leq K_5\ e^{-\gamma t}\limsup_{n\rightarrow\infty}\left(
\left\vert u_{0}^{n}-u_{0}\right\vert ^{2}+\left\Vert
\phi^{n}-\phi\right\Vert _{L_{V}^{2}}^{2}\right) ~~ \text{ for all }~~t\geq0,
\label{Conv3}
\end{equation}
where $K_5=\frac{1}{m}((\gamma+\delta)^2+1)$. Moreover, if $\left(
u_{0}^{n},\phi^{n}\right) \rightarrow\left( u_{0},\phi\right) $ strongly in $%
X$ as $n\rightarrow \infty$, then%
\begin{equation}
u^{n}(\cdot)\rightarrow u(\cdot)~~\text{ in }~~ L^{2}\left( 0,T;V \right) ~~%
\text{ for all }~~ T>0;  \label{Conv5}
\end{equation}%
\begin{equation}
u_t^{n}(\cdot)\rightarrow u_t(\cdot)~~ \text{ in }~~L_{V}^{2}~~ \text{ for
all }~~t\geq 0.  \label{Conv6}
\end{equation}
\end{lemma}

\noindent\textbf{Proof.} Let $T>0$ be arbitrary. In view of (\ref{IneqX})
and integrating in (\ref{Ineq0}) over $\left(0,T\right)$, we deduce that $%
u^{n}$ is bounded in $L^{\infty}(0,T;H)$, $L^{2}(0, T;V)$ and $L^{2p}(0,
T;L^{2p}\left( \Omega\right) )$, $\eta_{n}^{t}$ is bounded in $L^{\infty}(0,
T;L_{\mu}^{2}\left( \mathbb{R}^{+};V\right) )$. Hence, passing to a
subsequence, we have 
\begin{align}
u^{n} & \rightarrow u~~\text{ weak-star in }~~L^{\infty}(0, T;H);
\label{Conv7} \\[1.0ex]
u^{n}& \rightarrow u~~\text{ weakly in }~~L^{2}(0, T;V );  \notag \\[1.0ex]
u^{n}& \rightarrow u~~\text{ weakly in }~~L^{2p}(0, T;L^{2p}\left(
\Omega\right) );  \notag \\[1.0ex]
\eta_{n}^{t} & \rightarrow\eta^{t}~~\text{ weak-star in }~~L^{\infty}(0,
T;L_{\mu }^{2}( \mathbb{R}^{+};V));  \notag
\end{align}
thus (\ref{Conv2}) holds. Also, by the same arguments in the proof of
Theorem~\ref{Existence} (see Appendix), we deduce 
\begin{align}
\frac{du^{n}}{dt} & \rightarrow\frac{du}{dt}~~\text{ weakly in }~~L^{2}( 0,
T; V^*) +L^{q}\left(0, T;L^{q}\left( \Omega\right) \right) ,  \label{Conv8}
\\[1.0ex]
f\left( u^{n}\right) & \rightarrow\chi~~\text{ weakly in }~~L^{q}\left( 0,
T;L^{q}\left( \Omega\right) \right) .  \notag
\end{align}
In view of (\ref{Conv2}) and (\ref{Conv8}), making use of the Compactness
Theorem \cite{Robinson} we infer that (\ref{Conv2B}) holds true. Thus, $%
u^{n}\left( t,x\right) \rightarrow u\left( t,x\right) $, $f\left(
u^{n}\left( t,x\right) \right) \rightarrow f\left( u\left( t,x\right)
\right) $ for a.a. $(t,x)\in (0, T)\times \Omega $, so Lemma 1.3 in \cite%
{Lions} implies that $\chi=f\left( u\right) .$

By proceeding as in the proof of Theorem~\ref{Existence}, we obtain that $%
z(\cdot)=(u(\cdot),\eta^{\cdot})$ is a solution to problem \eqref{eq3.3}
with initial value $z\left(0\right) =\left( u_{0},\mathcal{J}\phi\right)$.
Thanks to the uniqueness of solution, a standard argument implies that the
above convergences are true for the whole sequence.

Further, we will prove (\ref{Conv1}). In a standard way, we obtain 
\begin{equation*}
u_{n}\rightarrow u~~\text{ in }C([0,T],V^{\ast }+L^{q}(\mathcal{O})).
\end{equation*}%
Then, if $t_{n}\rightarrow t_{0}$, $t_{n}\in \lbrack 0,T],\ t_{0}\in \lbrack
0,T]$, we infer%
\begin{equation*}
u_{n}(t_{n})\rightarrow u(t_{0})\text{ weakly in }H,
\end{equation*}%
and 
\begin{equation}
\left\vert u(t_{0})\right\vert \leq \liminf_{n\rightarrow \infty }\
\left\vert u_{n}(t_{n}))\right\vert .  \label{LimInf}
\end{equation}%
We need to prove that $u^{n}(t_{n})\rightarrow u(t_{0})$ strongly in $H$
when $t_{0}>0$. By Corollary \ref{CorExistSol}, we know that $u^{n}$ are
weak solutions of the equation%
\begin{equation*}
u_{t}^{n}-a(l(u))\Delta u^{n}-\int_{-\infty }^{t}k(t-s)\Delta
u^{n}ds+f(u^{n})=h.
\end{equation*}%
Then,%
\begin{equation*}
\frac{1}{2}\frac{d}{dt}|u^{n}(t)|^{2}+m\Vert u^{n}\Vert
^{2}+(f(u^{n}),u^{n})\leq \left( \int_{-\infty }^{t}k(t-s)\Delta
u^{n}(s)ds,u^{n}(t)\right) +(h,u^{n}(t)).
\end{equation*}%
By (\ref{fCondition}) and the Young inequality, we obtain 
\begin{equation*}
\frac{d}{dt}|u^{n}|^{2}+m\Vert u^{n}\Vert ^{2}+f_{0}\Vert u^{n}\Vert
_{2p}^{2p}\leq 2a_{0}|\mathcal{O}|+\frac{1}{m\lambda _{1}}%
|h|^{2}+2\int_{-\infty }^{t}k(t-s)\Vert u^{n}(s)\Vert ds\Vert u^{n}(t)\Vert .
\end{equation*}%
On the other hand,%
\begin{equation*}
\int_{-\infty }^{t}k(t-s)\Vert u^{n}(s)\Vert ds=\int_{-\infty
}^{0}k(t-s)\Vert u^{n}(s)\Vert ds+\int_{0}^{t}k(t-s)\Vert u^{n}(s)\Vert
ds=I_{1}+I_{2}.
\end{equation*}%
From $\gamma <\delta $ and Remark \ref{rem3.2}, $I_{1}$ is estimated by 
\begin{equation*}
\begin{split}
I_{1}& =\int_{-\infty }^{0}k(t-s)e^{-\frac{\gamma s}{2}}e^{\frac{\gamma s}{2}%
}\Vert \phi ^{n}(s)\Vert ds\leq \left( \int_{-\infty
}^{0}k^{2}(t-s)e^{-\gamma s}ds\right) ^{\frac{1}{2}} \\[0.8ex]
& ~~\times \left( \int_{-\infty }^{0}e^{\gamma s}\Vert \phi ^{n}\Vert
^{2}ds\right) ^{\frac{1}{2}}\leq \Vert \phi ^{n}\Vert _{L_{V}^{2}}\left(
\int_{t}^{\infty }k(s)e^{-\gamma (t-s)}ds\right) ^{\frac{1}{2}}M_{1}^{\frac{1%
}{2}} \\[0.8ex]
& \leq \frac{\Vert \phi ^{n}\Vert _{L_{V}^{2}}M_{1}^{\frac{1}{2}}}{\delta ^{%
\frac{1}{2}}}\left( \int_{t}^{\infty }\mu (s)e^{-\gamma (t-s)}ds\right) ^{%
\frac{1}{2}} \\[0.8ex]
& \leq \frac{\left\Vert \phi ^{n}\right\Vert _{L_{V}^{2}}M_{1}^{\frac{1}{2}}%
}{\delta ^{\frac{1}{2}}}\left( \int_{t}^{\infty }\mu (t)e^{-\delta
(s-t)}e^{-\gamma (t-s)}ds\right) ^{\frac{1}{2}} \\[0.8ex]
& \leq \frac{M_{1}^{\frac{1}{2}}\mu ^{\frac{1}{2}}(t)\Vert \phi ^{n}\Vert
_{L_{V}^{2}}}{\delta ^{\frac{1}{2}}(\delta -\gamma )^{\frac{1}{2}}}.
\end{split}%
\end{equation*}%
For $I_{2}$, by means of the property $k(t)\leq M_{1}$ and the boundedness
of $u^{n}$ in $L^{2}(0,T;V)$, there exists a constant $M_{2}$ such that 
\begin{equation*}
I_{2}\leq M_{1}\int_{0}^{t}\Vert u^{n}(s)\Vert ds\leq M_{2}\sqrt{t}.
\end{equation*}%
Thus,%
\begin{align*}
& \frac{d}{dt}|u^{n}(t)|^{2}+m\Vert u^{n}(t)\Vert ^{2}+f_{0}\Vert
u^{n}(t)\Vert _{2p}^{2p} \\
& \leq 2a_{0}|\mathcal{O}|+\frac{1}{m\lambda _{1}}|h|^{2}+2\left( \frac{%
M_{1}^{\frac{1}{2}}\mu ^{\frac{1}{2}}(t)\Vert \phi ^{n}\Vert _{L}{}_{V}^{2}}{%
\delta ^{\frac{1}{2}}(\delta -\gamma )^{\frac{1}{2}}}+M_{2}\sqrt{T}\right)
\Vert u^{n}(t)\Vert .
\end{align*}%
Therefore,%
\begin{align*}
& \frac{d}{dt}|u^{n}(t)|^{2}+\frac{m}{2}\Vert u^{n}(t)\Vert ^{2}+f_{0}\Vert
u^{n}(t)\Vert _{2p}^{2p} \\
& \leq 2a_{0}|\mathcal{O}|+\frac{1}{m\lambda _{1}}|h|^{2}+\frac{4M_{1}\mu
(t)\Vert \phi ^{n}\Vert _{L_{V}^{2}}^{2}}{\delta m(\delta -\gamma )}+\frac{%
4(M_{2})^{2}T}{m}.
\end{align*}%
The function $u$ satifies the same inequality but replacing $\phi ^{n}$ by $%
\phi $. We define the functions 
\begin{equation*}
J_{n}(t)=|u^{n}(t)|^{2}-2a_{0}|\mathcal{O}|t-\frac{4(M_{2})^{2}T}{m}t-\frac{%
\left\vert h\right\vert ^{2}}{m\lambda _{1}}t-\frac{4M_{1}\Vert \phi
^{n}\Vert _{L_{V}^{2}}^{2}}{\delta m(\delta -\gamma )}\int_{0}^{t}\mu (r)dr,
\end{equation*}%
\begin{equation*}
J(t)=|u(t)|^{2}-2a_{0}|\mathcal{O}|t-\frac{4(M_{2})^{2}T}{m}t-\frac{%
\left\vert h\right\vert ^{2}}{m\lambda _{1}}t-\frac{4M_{1}\Vert \phi \Vert
_{L_{V}^{2}}^{2}}{\delta m(\delta -\gamma )}\int_{0}^{t}\mu (r)dr.
\end{equation*}%
These functions are continuous and non-increasing on $[0,T]$, and 
\begin{equation*}
J_{n}(s)\rightarrow J(s)~~\text{ for a.a. }s\in \lbrack 0,T]~\text{as}%
~n\rightarrow \infty .
\end{equation*}%
Hence, there exists a sequence $\{t_{k}\}\in (0,t_{0})$ such that $%
t_{k}\rightarrow t_{0},$ when $k\rightarrow \infty $, and 
\begin{equation*}
\lim_{n\rightarrow \infty }J_{n}(t_{k})=J(t_{k}),\qquad \forall k\geq 1.
\end{equation*}%
Fix an arbitrary value $\epsilon >0$. From the continuity of $J$ on $[0,T]$,
there exists $k(\epsilon )\geq 1$ such that 
\begin{equation*}
|J(t_{k(\epsilon )})-J(t_{0})|\leq \frac{\epsilon }{2}.
\end{equation*}%
Now consider $n(\epsilon )\geq 1$ such that 
\begin{equation*}
t_{n}\geq t_{k(\epsilon )}~~\text{and}~~|J_{n}(t_{k(\epsilon
)})-J(t_{k(\epsilon )})|\leq \frac{\epsilon }{2},\quad \forall n\geq
n(\epsilon ).
\end{equation*}%
Then, since all $J_{n}$ are non-increasing, we deduce that 
\begin{equation*}
\begin{split}
J_{n}(t_{n})-J(t_{0})& \leq J_{n}(t_{k(\epsilon )})-J(t_{0})\leq
|J_{n}(t_{k(\epsilon )})-J(t_{0})| \\
& \leq |J_{n}(t_{k(\epsilon )})-J(t_{k(\epsilon )})|+|J(t_{k(\epsilon
)})-J(t_{0})|\leq \epsilon ,\qquad \forall n\geq n(\epsilon ).
\end{split}%
\end{equation*}%
As $\epsilon >0$ is arbitrary, we obtain 
\begin{equation*}
\limsup_{n\rightarrow \infty }J_{n}(t_{n})\leq J(t_{0}).
\end{equation*}%
Thus, 
\begin{equation}
\limsup_{n\rightarrow \infty }\ |u^{n}(t_{n})|\leq |u(t_{0})|.  \label{eq911}
\end{equation}%
Therefore, (\ref{LimInf}) and (\ref{eq911}) imply that $u^{n}(s_{n})%
\rightarrow u(t_{0})$ strongly in $H$, and (\ref{Conv1}) holds true.

Define the functions $w^{n}=z^{n}-z,\ \beta _{n}^{t}=\eta _{n}^{t}-\eta ^{t}$%
, similarly to the uniqueness step in the proof of Theorem~\ref{Existence},
Step 5 in Appendix, we have 
\begin{align}
& \frac{d}{dt}\left\Vert w^{n}\right\Vert _{\mathcal{H}}^{2}+2(((\beta
_{n}^{t})^{\prime },\beta _{n}^{t}))_{\mu }  \label{wn} \\
& \leq -2\int_{\Omega }\left( f\left( u^{n}\right) -f\left( u\right) \right)
\left( u^{n}-u\right) dx-\int_{\Omega }(a\left( l\left( u^{n}\right) \right)
\nabla u^{n}-a(l\left( u\right) )\nabla u)\text{\textperiodcentered }\nabla
\left( u^{n}-u\right) dx.  \notag
\end{align}%
Since $a$ is locally Lipschitz, by \eqref{eq2.2} and the Young inequality,
we have%
\begin{align}
& -2\int_{\Omega }(a\left( l\left( u^{n}\right) \right) \nabla
u^{n}-a(l\left( u\right) )\nabla u)\text{\textperiodcentered }\nabla \left(
u^{n}-u\right) dx  \notag \\
& =-2\int_{\Omega }a\left( l\left( u^{n}\right) \right) \left\vert \nabla
\left( u^{n}-u\right) \right\vert ^{2}dx-2\left( a\left( l\left(
u^{n}\right) \right) -a\left( l\left( u\right) \right) \right) \int_{\Omega
}\nabla u\text{\textperiodcentered }\nabla \left( u^{n}-u\right) dx  \notag
\\
& \leq -2m\left\Vert u^{n}-u\right\Vert ^{2}+2L_{a}\left( R\right)
\left\vert l\right\vert \left\vert u^{n}-u\right\vert \left\Vert
u\right\Vert \left\Vert u^{n}-u\right\Vert  \notag \\
& \leq \left( \alpha -2m\right) \left\Vert u^{n}-u\right\Vert ^{2}+\frac{%
L_{a}^{2}\left( R\right) \left\vert l\right\vert ^{2}}{\alpha }\left\vert
u^{n}-u\right\vert ^{2}\left\Vert u\right\Vert ^{2},  \label{Ineqa}
\end{align}%
where $\alpha \leq \left( m\lambda _{1}-\gamma \right) /\lambda _{1}$, and
for all $n\geq 1$, $t\geq 0$, we choose $R>0$ such that $\left\vert
u^{n}\left( t\right) \right\vert ,\ \left\vert u\left( t\right) \right\vert
\leq R$ (cf. \eqref{Conv1}). By (\ref{IneqEta}), (\ref{wn}) and (\ref{Ineqa}%
), we have%
\begin{align*}
& \frac{d}{dt}\left\Vert w^{n}\right\Vert _{\mathcal{H}}^{2}+\gamma
\left\Vert w^{n}\right\Vert _{\mathcal{H}}^{2}+m\left\Vert
u^{n}-u\right\Vert ^{2} \\
& \leq \frac{d}{dt}\left\Vert w^{n}\right\Vert _{\mathcal{H}}^{2}+\left(
2m-\alpha \right) \left\Vert u^{n}-u\right\Vert ^{2}+\delta \int_{0}^{\infty
}\mu (s)|\nabla \beta _{n}^{t}(s)|^{2}ds \\
& \leq \frac{L_{a}^{2}\left( R\right) \left\vert l\right\vert ^{2}}{\alpha }%
\left\vert u^{n}-u\right\vert ^{2}\left\Vert u\right\Vert ^{2}-2\int_{\Omega
}\left( f\left( u^{n}\right) -f\left( u\right) \right) \left( u^{n}-u\right)
dx,
\end{align*}%
where we have used that {$\gamma \leq \min \{\left( m-\alpha \right) \lambda
_{1},\delta \}$} by the choice of $\alpha $. Multiplying by $e^{\gamma t}$
on both sides of the above inequality and integrating over $\left(
0,t\right) $, we obtain%
\begin{align*}
& \left\Vert w^{n}\left( t\right) \right\Vert _{\mathcal{H}%
}^{2}+m\int_{0}^{t}e^{-\gamma (t-s)}\left\Vert u^{n}-u\right\Vert ^{2}ds \\
& \leq e^{-\gamma t}\left\Vert w^{n}\left( 0\right) \right\Vert _{\mathcal{H}%
}^{2}+\frac{L_{a}^{2}\left( R\right) \left\vert l\right\vert ^{2}}{\alpha }%
\int_{0}^{t}e^{-\gamma (t-s)}\left\vert u^{n}-u\right\vert ^{2}\left\Vert
u\right\Vert ^{2}ds \\
& \quad -2\int_{0}^{t}e^{-\gamma (t-s)}\int_{\Omega }\left( f\left(
u^{n}\right) -f\left( u\right) \right) \left( u^{n}-u\right) dxds.
\end{align*}%
On the one hand, by (\ref{Conv1}), we know that $\left\vert u^{n}\left(
s\right) -u\left( s\right) \right\vert ^{2}\left\Vert u\left( s\right)
\right\Vert ^{2}\rightarrow 0$ for a.e. $s\in \left( 0,t\right) $. On the
other hand, $e^{-\gamma (t-s)}\left\vert u^{n}\left( s\right) -u\left(
s\right) \right\vert ^{2}\left\Vert u\left( s\right) \right\Vert ^{2}$ is
bounded by the integrable function $4R^{2}e^{-\gamma (t-s)}\left\Vert
u\left( s\right) \right\Vert ^{2}$. Hence, Lebesgue's theorem implies that 
\begin{equation*}
\int_{0}^{t}e^{-\gamma (t-s)}\left\vert u^{n}-u\right\vert ^{2}\left\Vert
u\right\Vert ^{2}ds\rightarrow 0~~\text{ as }~~n\rightarrow \infty .
\end{equation*}%
Since $f\left( u^{n}\right) \rightarrow f\left( u\right) $ weakly in $%
L^{q}\left( 0,T;L^{q}\left( \Omega \right) \right) $, it follows that%
\begin{equation*}
\int_{0}^{t}e^{-\gamma (t-s)}\int_{\Omega }\left( f\left( u^{n}\right)
-f\left( u\right) \right) udxds\rightarrow 0~~\mbox{as}~~n\rightarrow \infty
.
\end{equation*}%
Furthermore, as $f\left( u^{n}\left( t,x\right) \right) u^{n}\left(
t,x\right) \geq -\kappa _{1}+\kappa _{2}|u^{n}(t,x)|^{2p}$ (see %
\eqref{eq3.11}) and $u^{n}\left( t,x\right) \rightarrow u\left( t,x\right) $%
, $f\left( u^{n}\left( t,x\right) \right) \rightarrow f\left( u\left(
t,x\right) \right) $ for a.a. $\left( t,x\right) \in (0,T]\times \Omega $,
Lebesgue-Fatous's theorem implies%
\begin{align*}
& \limsup_{n\rightarrow \infty }\left( -2\int_{0}^{t}e^{-\gamma
(t-s)}\int_{\Omega }f\left( u^{n}\right) u^{n}dxds\right) \\
& \leq -2\liminf_{n\rightarrow \infty }\int_{0}^{t}e^{-\gamma
(t-s)}\int_{\Omega }f\left( u^{n}\right) u^{n}dxds \\
& \leq -2\int_{0}^{t}e^{-\gamma (t-s)}\int_{\Omega }\liminf_{n\rightarrow
\infty }f\left( u^{n}\right) u^{n}dxds \\
& =2\int_{0}^{t}e^{-\gamma (t-s)}\int_{\Omega }f\left( u\right) udxds.
\end{align*}%
This inequality, together with 
\begin{equation}
\int_{0}^{t}e^{-\gamma (t-s)}\int_{\Omega }f\left( u\right) \left(
u^{n}-u\right) dxds\rightarrow 0~~\mbox{as}~~n\rightarrow \infty ,
\label{nimo}
\end{equation}%
shows that%
\begin{equation*}
\limsup_{n\rightarrow \infty }\left( -2\int_{0}^{t}e^{-\gamma
(t-s)}\int_{\Omega }(f\left( u^{n}\right) -f(u))u^{n}dxds\right) \leq 0~~%
\mbox{as}~~n\rightarrow \infty .
\end{equation*}%
Notice that \eqref{nimo} follows from the facts $f\left( u(\cdot )\right)
\in L^{q}\left( 0,T;L^{q}\left( \Omega \right) \right) $ and $%
u^{n}\rightarrow u$ weakly in $L^{2p}\left( 0,T;L^{2p}\left( \Omega \right)
\right) $.

Collecting all inequalities and using \eqref{Kmu}, 
\begin{align*}
& \limsup_{n\rightarrow\infty}\int_{0}^{t}e^{-\gamma(t-s)}\left\Vert
u^{n}(s)-u(s)\right\Vert ^{2}ds \\
& \leq\frac{1}{m}e^{-\gamma t}\limsup_{n\rightarrow\infty}\left\Vert
w^{n}\left(0\right) \right\Vert _{\mathcal{H}}^{2} \\
& =\frac{1}{m}e^{-\gamma t}\limsup_{n\rightarrow\infty}\left( \left\vert
u^{n}\left(0\right) -u_{0}\right\vert ^{2}+\int_{0}^{\infty}\mu\left(
s\right) \|\beta_{n}^{0}\left( s\right) \| ^{2}ds\right) \\
& \leq\frac{1}{m}e^{-\gamma t}\limsup_{n\rightarrow\infty}\left( \left\vert
u^{n}\left(0\right) -u_{0}\right\vert ^{2}+K_\mu\int_{-\infty}^{0}e^{\gamma
s}\left\Vert \phi^{n}\left( s\right) -\phi\left( s\right) \right\Vert
^{2}ds\right).
\end{align*}
Finally, (\ref{Conv3}) follows from%
\begin{align*}
\left\Vert u_{t}^{n}-u_{t}\right\Vert _{L_{V}^{2}}^{2} &
=\int_{-t}^{0}e^{\gamma s}\left\Vert u^{n}\left( t+s\right) -u\left(
t+s\right) \right\Vert ^{2}ds+\int_{-\infty}^{-t}e^{\gamma s}\left\Vert
u^{n}\left( t+s\right) -u\left( t+s\right) \right\Vert ^{2}ds \\
& =\int_{0}^{t}e^{-\gamma(t-s)}\left\Vert u^{n}\left( s\right) -u\left(
s\right) \right\Vert ^{2}ds+e^{-\gamma t}\int_{-\infty}^{0}e^{\gamma
s}\left\Vert \phi^{n}\left( s\right) -\phi\left( s\right) \right\Vert ^{2}ds.
\end{align*}
If $\left( u_{0}^{n},\phi^{n}\right) \rightarrow\left( u_{0},\phi\right) $
in $X$, then (\ref{Conv3}) implies (\ref{Conv5}) and (\ref{Conv6}). $\Box$

\medskip

\begin{remark}
The convergence (\ref{Conv1}) was stated in Lemma 3.9 from \cite{X2}.
However, the proof of this result is incorrect there and we provided here a
correct one.
\end{remark}

As a consequence, we obtain the continuous dependence with respect to the
initial data.

\begin{corollary}
\label{Continuity}Assume conditions of Theorem \ref{Existence} are true.
Then, for any $t\geq0$, the mapping $\left( u_{0},\phi\right) \mapsto
S(t)\left( u_{0},\phi\right) $ is continuous.
\end{corollary}

Finally, we are ready to prove the asymptotic compactness of the semigroup.

\begin{lemma}
\label{AC} Under assumptions of Theorem~\ref{Existence}, the semigroup $S$
is asymptotically compact.
\end{lemma}

\noindent \textbf{Proof.} Let $B\subset X$ be a bounded set, we need to
prove that for any sequences $\{(y_{n},\phi _{n})\}_{n\in \mathbb{N}}\subset
B$ and $t_{n}\rightarrow +\infty $ as $n\rightarrow +\infty $, the sequence $%
\{S(t_{n})(y_{n},\phi _{n})\}_{n\in \mathbb{N}}$ is relatively compact.
Recall that 
\begin{equation*}
S(t_{n})(y_{n},\phi _{n})=(u(t_{n};0,(y_{n},\mathcal{J}\phi
_{n})),u_{t_{n}}(\cdot ;0,(y_{n},\mathcal{J}\phi
_{n}))):=(u^{n}(t_{n}),u_{t_{n}}^{n}(\cdot ))
\end{equation*}%
Pick now $T>0$, and assume that $t_{n}>T$ for all $n\in \mathbb{N}$ (there
is no loss of generality in assuming this since $t_{n}\rightarrow +\infty $%
). Now we can define $v^{n}(t)=u^{n}\left( t+t_{n}-T\right) $, observe that $%
v^{n}(T)=u^{n}(t_{n})$ and $%
v_{T}^{n}(t)=v^{n}(T+t)=u^{n}(t+t_{n})=u_{t_{n}}^{n}(t).$ Therefore 
\begin{equation*}
S(t_{n})(y_{n},\phi _{n})=(u^{n}(t_{n}),u_{t_{n}}^{n}(\cdot
))=(v^{n}(T),v_{T}^{n}(\cdot )).
\end{equation*}%
Let us denote now 
\begin{equation*}
\mathcal{Y}_{n}=(v^{n}(T),v_{T}^{n})=(u^{n}(t_{n}),u_{t_{n}}^{n}(\cdot )),\
\xi _{n}^{T}=(v^{n}(0),v_{0}^{n}(\cdot
))=(u^{n}(t_{n}-T),u_{t_{n}-T}^{n}(\cdot )).
\end{equation*}%
By Lemma \ref{Est1}, the sequences $\{\mathcal{Y}_{n}\}$,\ $\{\xi _{n}^{T}\}$
are bounded in $X$, so up to a subsequence $\mathcal{Y}_{n}\rightarrow 
\mathcal{Y}:=(y,\phi )$, $\xi _{n}^{T}\rightarrow \xi ^{T}$ weakly in $X$.
In addition, by Lemma \ref{ConvergSol}, $\mathcal{V}(t):=S(t)\xi
^{T}=(v(t),v_{t}(\cdot ))$ satisfies (\ref{Conv1})-(\ref{Conv3}). It follows
from the above convergences that, $\phi =v_{T}$ in $L_{V}^{2}$ and $%
y=v_{T}\left( 0\right) ,\ \phi \left( s\right) =v_{T}\left( s\right) $ for
almost all $s\in (-\infty ,0)$. Also, in view of (\ref{Conv1}) we infer that 
\begin{equation*}
u^{n}(t_{n})=v^{n}(T)\rightarrow v\left( T\right) =y.
\end{equation*}%
Hence, in order to prove that $\mathcal{Y}_{n}\rightarrow \mathcal{Y}$ in $X$%
, it remains to check that $u_{t_{n}}^{n}(\cdot )\rightarrow \phi $ in $%
L_{V}^{2}$ (up to a subsequence). Notice that $u_{t_{n}}^{n}(\cdot
)=v_{T}^{n}$ for all $t_{n}>T$ and $v_{T}=\phi $. Thanks to (\ref{Conv3}) we
have, for each $T>0$,%
\begin{align*}
\limsup_{n\rightarrow \infty }\ \left\Vert u_{t_{n}}^{n}(\cdot )-\phi
\right\Vert _{L_{V}^{2}}^{2}& =\limsup_{n\rightarrow \infty }\ \left\Vert
v_{T}^{n}-v_{T}\right\Vert _{L_{V}^{2}}^{2} \\
& \leq K_{5}\ e^{-\gamma (T-\tau )}\limsup_{n\rightarrow \infty }\left(
\left\Vert \xi _{n}^{T}-\xi ^{T}\right\Vert _{X}^{2}\right) \\
& \leq \widetilde{K}e^{-\gamma T},
\end{align*}%
where the last inequality follows from Lemma \ref{Est1}. For every $k>0$,
there exists $T:=T\left( k\right) $ such that for all $T\geq T(k)$, 
\begin{equation*}
\limsup_{n\rightarrow \infty }\ \left\Vert u_{t_{n}}^{n}(\cdot )-\phi
\right\Vert _{L_{V}^{2}}^{2}=\limsup_{n\rightarrow \infty }\ \left\Vert
v_{T}^{n}-v_{T}\right\Vert _{L_{V}^{2}}^{2}\leq \frac{1}{k}.
\end{equation*}%
Taking $k\rightarrow \infty $ and using a diagonal argument, we obtain that
there exists a subsequence $\{u_{t_{n_{k}}}^{n_{k}}(\cdot )\}$ such that $%
u_{t_{n_{k}}}^{n_{k}}(\cdot )\rightarrow \phi $ in $L_{V}^{2}$. $\Box $

\bigskip

By Corollaries \ref{Absorbing}, \ref{Continuity} and Lemma \ref{AC} the
general theory of attractors (see \cite[Theorem 3.1]{Lad}) implies the
following result.

\begin{theorem}
Under the assumptions of Theorem~\ref{Existence}, the semigroup $S$
possesses a global connected attractor $\mathcal{A}\subset X.$
\end{theorem}

As a straightforward consequence of the previous results, we can provide
information for the local problem analyzed, amongst others, in the papers 
\cite{G, P2, C} by simply assuming that $a(\cdot)$ is a constant function.

\begin{corollary}
Under the hypotheses of Theorem~\ref{Existence}, assume also that $%
a(t)=k_0>0 $ for all $t\geq0$. Then the local problem \eqref{eq0} poseesses
a global connected attractor $\mathcal{A}\subset X.$
\end{corollary}

\section*{Acknowledgements}

The research has been partially supported by the Spanish Ministerio de
Ciencia, Innovaci\'on y Universidades (MCIU), Agencia Estatal de
Investigaci\'on (AEI) and Fondo Europeo de Desarrollo Regional (FEDER) under
the project PGC2018-096540-B-I00, and by Junta de Andaluc\'{\i}a (Consejer%
\'{\i}a de Econom\'{\i}a y Conocimiento) and FEDER under projects US-1254251
and P18-FR-4509.\newline

\section*{Appendix}

\label{appendix}

\textbf{Proof of Theorem~\ref{Existence}.} We follow a standard
Faedo-Galerkin method. Recall that there exists a smooth orthonormal basis $%
\{w_j\}_{j=1}^{\infty}$ in $H$ which also belongs to $V\cap L^{2p}(\Omega)$ (%
\cite[Proposition 2.3]{C}). Let us take a complete set of normalized
eigenfunctions for $-\Delta$ in $V$, such that $-\Delta w_j=\lambda_j w_j$ ($%
\lambda_j$ the eigenvalue corresponding to $w_j$). Next we will select an
orthonormal basis $\{\zeta_j\}_{j=1}^{\infty} $ of $L^2_{\mu}(\mathbb{R}%
^+;V) $ which also belongs to $\mathcal{D}(\mathbb{R}^+;V)$.

The proof is divided into $6$ steps.

\textbf{Step 1.} (Faedo-Galerkin scheme) Fix $T>\tau$, for a given integer $%
n $, denote by $P_n$ and $Q_n$ the projections on the subspaces 
\begin{equation*}
\mbox{span}\{w_1,\cdots,w_n\}\subset V\qquad \mbox{and}\qquad \mbox{span}%
\{\zeta_1,\cdots,\zeta_n\}\subset L^2_{\mu}(\mathbb{R}^+;V),
\end{equation*}
respectively. We look for a function $z_n=(u_n,\eta^t_n)$ of the form 
\begin{equation*}
u_n(t)=\sum_{j=1}^nb_j(t)w_j \qquad \mbox{and}\qquad
\eta^t_n(s)=\sum_{j=1}^nc_j(t)\zeta_j(s),
\end{equation*}
satisfying 
\begin{alignat}{2}  \label{eq3.6}
\begin{cases}
(\partial_tz_n,(w_k,\zeta_j))_{\mathcal{H}}=(\mathcal{L}z_n,(w_k,\zeta_j))+(%
\mathcal{G}(z),(w_k,\zeta_j)), \\ 
z_n|_{t=\tau}=(P_nu_0,Q_n\eta_0),%
\end{cases}
\begin{aligned} & k,j=0,\cdots,n,\\ &\qquad \end{aligned} & & 
\end{alignat}
for a.e. $\tau\leq t\leq T$, where $w_0$ and $\zeta_0$ are the zero vectors
in the respective spaces. Taking $(w_k,\zeta_0)$ and $(w_0,\zeta_k)$ in (\ref%
{eq3.6}), applying the divergence theorem, we derive a system of ODE in the
variables 
\begin{alignat}{2}  \label{eq3.7}
\begin{cases}
\frac{d}{dt}b_k(t)=-\lambda_ka(l(\sum_{j=1}^n{b_j(t)w_j}))b_k-%
\sum_{j=1}^nc_j((\zeta_j,w_k))_{\mu}-({f(\sum_{j=1}^nb_j(t)w_j)}%
,w_k)+(g,w_k), \\[1.0ex] 
\frac{d}{dt}c_k(t)=\sum_{j=1}^nb_j((w_j,\zeta_k))_{\mu}-\sum_{j=1}^nc_j((%
\zeta_j^{\prime},\zeta_k))_{\mu}, \\[1.0ex] 
\end{cases}
& & 
\end{alignat}
subject to the initial conditions, 
\begin{equation}  \label{eq3.8}
b_k(\tau)=(u_0,w_k),\qquad c_k(\tau)=((\eta_0,\zeta_k))_{\mu}.
\end{equation}
According to the standard existence theory for ODE, there exists a
continuous solution of (\ref{eq3.7})-(\ref{eq3.8}) on some interval $%
(\tau,t_n)$. Then a priori estimates imply $t_n=\infty$.

\textbf{Step 2.} (Energy estimate) Multiplying the first equation of (\ref%
{eq3.7}) by $b_{k}$ and the second one by $c_{k}$, summing over $k$ {($%
k=1,2,\cdots ,n$)} and adding the results, we have 
\begin{equation}
\frac{1}{2}\frac{d}{dt}\Vert z_{n}\Vert _{\mathcal{H}}^{2}=(\mathcal{L}%
z_{n},z_{n})_{\mathcal{H}}+(\mathcal{G}(z_{n}),z_{n})_{\mathcal{H}}.
\label{eq3.9}
\end{equation}%
On the one hand, by the divergence theorem, 
\begin{equation*}
\left( \int_{0}^{\infty }\mu (s)\Delta \eta _{n}^{t}(s)ds,u_{n}\right)
=-\int_{0}^{\infty }\mu (s)\int_{\Omega }\nabla \eta _{n}^{t}(s)\cdot \nabla
u_{n}(s)dxds=-((u_{n},\eta _{n}^{t}))_{\mu },
\end{equation*}%
therefore, 
\begin{equation}
(\mathcal{L}z_{n},z_{n})_{\mathcal{H}}=-a(l(u_{n}))|\nabla
u_{n}|^{2}-(((\eta _{n}^{t})^{\prime },\eta _{n}^{t}))_{\mu }.
\label{eq3.10}
\end{equation}%
On the other hand, (\ref{eq3.4}) yields that there exists a constant $a_{0}$%
, such that 
\begin{equation}
f(u)\cdot u\geq \frac{1}{2}f_{0}u^{2p}-a_{0},  \label{fCondition}
\end{equation}%
hence, 
\begin{equation}
(\mathcal{G}(z_{n}),z_{n})_{\mathcal{H}}=(-f(u_{n})+g,u_{n})\leq -\frac{1}{2}%
f_{0}|u_{n}|_{2p}^{2p}+a_{0}|\Omega |+(g,u_{n}).  \label{eq3.11}
\end{equation}%
{It follows from} (\ref{eq2.2}), (\ref{eq3.9})-(\ref{eq3.11}) and the Young
inequality that 
\begin{equation}
\frac{d}{dt}\Vert z_{n}\Vert _{\mathcal{H}}^{2}+2m|\nabla
u_{n}|^{2}+2(((\eta _{n}^{t})^{\prime },\eta _{n}^{t}))_{\mu
}+f_{0}|u_{n}|_{2p}^{2p}\leq 2a_{0}|\Omega |+\frac{1}{m\lambda _{1}}%
|g|^{2}+m|\nabla u_{n}|^{2}.  \label{eq611-3}
\end{equation}%
Integration by parts and $(h_{1})$ yield that, 
\begin{equation*}
2(((\eta _{n}^{t})^{\prime },\eta _{n}^{t}))_{\mu }=-\int_{0}^{\infty }\mu
^{\prime }(s)|\nabla \eta _{n}^{t}(s)|^{2}ds\geq 0,
\end{equation*}%
thus the third term of the right hand side of \eqref{eq611-3} can be
neglected, we obtain 
\begin{equation*}
\frac{d}{dt}\Vert z_{n}\Vert _{\mathcal{H}}^{2}+m|\nabla
u_{n}|^{2}+f_{0}|u_{n}|_{2p}^{2p}\leq 2a_{0}|\Omega |+\frac{1}{m\lambda _{1}}%
|g|^{2}.
\end{equation*}%
%
%
%
%
Integrating the above inequality between $\tau $ and $t$, $t\in (\tau ,T]$,
we have 
\begin{equation}
\Vert z_{n}(t)\Vert _{\mathcal{H}}^{2}+\int_{\tau }^{t}\left[ m\Vert
u_{n}\Vert ^{2}+f_{0}|u_{n}|_{2p}^{2p}\right] dr\leq \Vert z_{0}\Vert _{%
\mathcal{H}}^{2}+\Lambda (T-\tau ),  \label{3.11}
\end{equation}%
where we have used the notation $\Lambda :=2a_{0}|\Omega |+\frac{1}{m\lambda
_{1}}|g|^{2}$. Therefore, it arrives that 
\begin{equation*}
\begin{split}
& u_{n}\quad \mbox{is bounded in}\quad L^{\infty }(\tau ,T;H)\cap L^{2}(\tau
,T;V)\cap L^{2p}(\tau ,T;L^{2p}(\Omega )), \\
& \eta _{n}\quad \mbox{is bounded in}\quad L^{\infty }(\tau ,T;L_{\mu }^{2}(%
\mathbb{R}^{+};V)).
\end{split}%
\end{equation*}%
Passing to a subsequence, there exists a function $z=(u,\eta )$ such that 
\begin{equation}
\left\{ 
\begin{array}{rcl}
\begin{aligned} &u_n\rightarrow u\qquad \mbox{weak-star in}\qquad
L^{\infty}(\tau,T;H);\\ &u_n\rightarrow u\qquad ~\mbox{weakly in}\qquad
\quad L^2(\tau,T;V);\\ &u_n\rightarrow u\qquad ~\mbox{weakly in}\qquad \quad
L^{2p}(\tau,T;L^{2p}(\Omega));\\ &\eta^t_n\rightarrow \eta^t \qquad
\mbox{weak-star in}\qquad L^{\infty}(\tau,T;L^2_{\mu}(\mathbb{R}^+;V)).
\end{aligned} &  & 
\end{array}%
\right.  \label{eq3.13}
\end{equation}

\textbf{Step 3.} (Passage to limit) For a fixed integer $m$, choose a
function 
\begin{equation*}
v=(\sigma, \pi)\in\mathcal{D}((\tau,T); V\cap L^{2p}(\Omega))\times \mathcal{%
D}((\tau,T);\mathcal{D}(\mathbb{R}^+;V))
\end{equation*}
of the form 
\begin{equation*}
\sigma(t)=\sum_{j=1}^{m}\tilde{b}_j(t)w_j\qquad \mbox{and}\qquad
\pi^t(s)=\sum_{j=1}^m\tilde{c}_j(t)\zeta_j(s),
\end{equation*}
where $\{\tilde{b}_j\}_{j=1}^m$ and $\{\tilde{c}_j\}_{j=1}^m$ are given
functions in $\mathcal{D}(\tau,T)$, then (\ref{eq3.6}) holds with $%
(\sigma,\pi)$ in place of $(\omega_k,\zeta_j)$.

Our main target is to prove problem (\ref{eq3.3}) has a solution in the weak
sense, i.e., for arbitrary $v\in \mathcal{D}((\tau,T); V\cap
L^{2p}(\Omega))\times \mathcal{D}((\tau,T);\mathcal{D}(\mathbb{R}^+;V))$
(here, specially, we pick up $v=(\sigma,\pi)\in \mathcal{D}(\tau,T)$ as a
test function), the following equality 
\begin{equation}  \label{eq3}
\begin{split}
\int_{\tau}^t(\partial_{r}z_n,v)_{\mathcal{H}}dr&=\int_{\tau}^t\bigg[%
-a(l(u_n))(\nabla u_n,\nabla\sigma)-((\eta_n^t,\sigma))_{\mu}-(f(u_n),\sigma)
\\
&\qquad\qquad +(g,\sigma)+((u_n,\pi^t))_{\mu}-\ll(\eta_n^t)^{\prime},\pi^t\gg%
\bigg]dr
\end{split}%
\end{equation}
holds in the sense of $\mathcal{D}^{\prime}(\tau,T)$. Here, we denote by $%
\ll\cdot,\cdot\gg$ the duality map between $H_{\mu}^1(\mathbb{R}^+;V)$ and
its dual space.

Firstly, using the same argument as in \cite[Theorem 2.7]{X1} and (\ref%
{eq3.13})$_2$, we know 
\begin{equation*}
\int_{\tau}^ta(l(u_n))(\nabla u_n,\nabla \sigma)dr \rightarrow
\int_{\tau}^ta(l(u))(\nabla u,\nabla \sigma)dr \quad \mbox{as}\quad
n\rightarrow \infty.
\end{equation*}
Similarly, by means of (\ref{eq3.13})$_4$ and (\ref{eq3.13})$_2$, we have 
\begin{equation*}
\int_{\tau}^t((\eta^t_n,\sigma))_{\mu} dr \rightarrow
\int_{\tau}^t((\eta^t,\sigma))_{\mu} dr \qquad \mbox{as}\quad n\rightarrow
\infty,
\end{equation*}
and 
\begin{equation*}
\int_{\tau}^t((u_n,\pi^t))_{\mu}dr\rightarrow
\int_{\tau}^t((u,\pi^t))_{\mu}dr \qquad \mbox{as}\quad n\rightarrow \infty,
\end{equation*}
respectively.

Secondly, we now show that 
\begin{equation*}
\lim_{n\rightarrow
\infty}\ll(\eta_n^t)^{\prime},\pi^t\gg=\ll(\eta^t)^{\prime},\pi^t\gg.
\end{equation*}
Notice that, for every $\upsilon\in L^2_{\mu}(\mathbb{R}^+;V)$, making use
of integration by parts, we derive 
\begin{equation}  \label{eq3.14}
\ll \upsilon^{\prime},\pi^t\gg=-\int_0^{\infty}\mu^{\prime}(s)(\nabla
\upsilon(s),\nabla\pi^t(s))ds-\int_0^{\infty}
\mu(s)(\nabla\upsilon(s),\nabla(\pi^t)^{\prime}(s))ds.
\end{equation}
Replacing $\upsilon$ by {$\eta^t_n$} in (\ref{eq3.14}), together with $(\ref%
{eq3.13})_4$, it is clear the right hand side of (\ref{eq3.14}) converges to 
$\ll (\eta^t)^{\prime},\pi^t\gg$ as $n\rightarrow \infty$.

Thirdly, we are going to prove that%
\begin{equation*}
f(u_{n})\rightarrow f(u)\text{ weakly in }L^{q}(\tau ,T;L^{q}(\Omega )).
\end{equation*}%
This follows from Lemma 1.3 in \cite[Chapter I]{Lions} if we show that 
\begin{equation*}
f(u_{n}(t,x))\rightarrow f(u(t,x))\qquad \mbox{for~ a.e.}~(t,x)\in (\tau
,T)\times \Omega ,
\end{equation*}%
and 
\begin{equation*}
|f(u_{n})|_{L^{q}((\tau ,T)\times \Omega )}\leq C,
\end{equation*}%
where $q=\frac{2p}{2p-1}\in (1,2)$, which is the conjugate exponent of $2p$
and the constant $C$ is independent of $n$. Observe that 
\begin{equation}
\begin{split}
\Vert \partial _{t}u_{n}\Vert _{L^{2}(\tau ,T;V^{\ast })+L^{q}(\tau
,T;L^{q}(\Omega ))}& \leq \Vert a(l(u_{n}))\Delta u_{n}\Vert _{L^{2}(\tau
,T;V^{\ast })}+\left\Vert \int_{0}^{\infty }\mu (s)\Delta \eta
_{n}^{t}(s)ds\right\Vert _{L^{2}(\tau ,T;V^{\ast })} \\[0.8ex]
& ~~+\Vert g\Vert _{V^{\ast }}+\Vert f(u)\Vert _{L^{q}(\tau ,T;L^{q}(\Omega
))}.
\end{split}
\label{eq3.15}
\end{equation}%
It follows from (\ref{eq3.4}), there exists a constant $K>0$ such that 
\begin{equation}
|f(u_{n})|^{q}\leq K(1+|u_{n}|^{2p}).  \label{eq3.16}
\end{equation}%
Together with (\ref{m}), (\ref{eq3.13}) and the assumption $g\in H$, we know
that $\{\partial _{t}u_{n}\}$ is bounded in $L^{2}(\tau ,T;V^{\ast
})+L^{q}(\tau ,T;L^{q}(\Omega ))$. Thus, up to a subsequence, we infer 
\begin{equation}
\partial _{t}u_{n}\rightarrow \tilde{u}\qquad \mbox{weakly in }\quad
L^{2}(\tau ,T;V^{\ast })+L^{q}(\tau ,T;L^{q}(\Omega )).  \label{eq3.17}
\end{equation}%
By a standard argument we infer that $\tilde{u}=u_{t}$. Since 
\begin{equation*}
L^{2}(\tau ,T;V^{\ast })+L^{q}(\tau ,T;L^{q}(\Omega ))\subset L^{q}(\tau
,T;V^{\ast }+L^{q}(\Omega ))
\end{equation*}%
and 
\begin{equation*}
L^{2}(\tau ,T;V)\subset L^{q}(\tau ,T;V),
\end{equation*}%
by (\ref{eq3.13}) and (\ref{eq3.17}), we deduce 
\begin{equation}
u_{n}\rightarrow u\qquad \mbox{weakly in}\quad W^{1,q}(\tau ,T;V^{\ast
}+L^{q}(\Omega ))\cap L^{q}(\tau ,T;V).  \label{eq3.18}
\end{equation}%
Applying a compactness argument \cite{Lions}, we derive that the injection 
\begin{equation*}
W^{1,q}(\tau ,T;V^{\ast }+L^{q}(\Omega ))\cap L^{q}(\tau
,T;V)\hookrightarrow L^{q}(\tau ,T;L^{q}(\Omega ))
\end{equation*}%
is compact. Therefore, (\ref{eq3.18}) implies that 
\begin{equation*}
u_{n}\rightarrow u\qquad \mbox{strongly in }\quad L^{q}(\tau ,T;L^{q}(\Omega
)).
\end{equation*}%
By the continuity of $f$ we obtain that (up to a subsequence) 
\begin{equation*}
f(u_{n}(t,x))\rightarrow f(u(t,x))\qquad \mbox{for~~a.e.}~(t,x)\in (\tau
,T)\times \Omega .
\end{equation*}%
In virtue of (\ref{eq3.16}), we obtain 
\begin{equation*}
|f(u_{n})|_{L^{q}((\tau ,T)\times \Omega )}^{q}=\int_{\tau }^{T}\int_{\Omega
}|f(u_{n})|^{q}dxdt\leq K|\Omega |(T-\tau )+K\int_{\tau
}^{T}|u_{n}|_{2p}^{2p}dt,
\end{equation*}%
which is bounded uniformly with respect to $n$.

Eventually, by a standard argument, we derive 
\begin{equation*}
\partial_t z_n\rightarrow z_t\quad \mbox{in}\quad \mathcal{D}%
^{\prime}(\tau,T;V\cap L^{2p}(\Omega))\times \mathcal{D}^{\prime}(\tau,T;%
\mathcal{D}(\mathbb{R}^+;V)).
\end{equation*}
Therefore, using a density argument, (\ref{eq3}) is proved by the previous
statements.

\textbf{Step 4.} (Continuity of solution) By (\ref{eq3.14})-(\ref{eq3.15}),
it is immediate to see that $z_t=(u_t,{\eta_t})$ fulfills 
\begin{equation*}
\begin{split}
& u_t\in L^2(\tau,T;V^*)+L^q(\tau,T;L^q(\Omega)); \\
&\eta_t\in L^2(\tau,T;H^{-1}_{\mu}(\mathbb{R}^+;V)),
\end{split}%
\end{equation*}
where $L^2(\tau,T;V^*)+L^q(\tau,T;L^q(\Omega))$ is the dual space of $%
L^2(\tau,T,V)\cap L^{2p}(\tau,T;L^{2p}(\Omega))$. Using a slightly modified
version of \cite[Lemma III.1.2]{T1}, together with (\ref{eq3.13}), we infer
that $u\in C([\tau,T];H)$.

As for the second component, by means of the same argument as \cite[Theorem,
Section 2]{C}, we obtain that {$\eta^t\in C([\tau,T];L_{\mu}^{2}(\mathbb{R}%
^+;V))$}. Thus, $z(\tau)$ makes sense, and the equality $z(\tau)=z_0$
follows from the fact that $(P_nu_0,Q_n\eta_0)$ converges to $z_0$ strongly.

\textbf{Step 5.} (Continuity with respect to the initial value and
uniqueness) Let $z_1=(u_1,\eta_1)$ and $z_2=(u_2,\eta_2)$ be the two
solutions of (\ref{eq3.5}) with initial data $z_{10}$ and $z_{20}$,
respectively. Due to the a priori estimates on the first component of
solutions $u$, see \eqref{3.11}, together with the fact that $u\in C(\tau,T;
H)$, we can ensure that there exists a bounded set $S\subset H$, such that $%
u_i(t)\in S$ for all $t\in[\tau, T]$ and $i=1,2$. In addition, taking into
account that $l\in \mathcal{L}(H;\mathbb{R})$, we have $\{l(u_i(t))\}_{t\in
[\tau, T]}\subset [-R,R]$ for $i=1,2$, for some $R>0$. Therefore, let $\bar{z%
}=z_1-z_2=(\bar{u},\bar{\eta})=(u_1-u_2,\eta_1-\eta_2)$ and $\bar{z}%
_0=z_{10}-z_{20}$. Thanks to \eqref{eq2.2}, the locally Lipschitz continuity
of function $a$ with Lipschitz constant $L_a(R)$ and the Poincar\'e
inequality, we have 
\begin{equation}  \label{eq3.20}
\begin{split}
\frac{d}{dt}\|\bar{z}\|^2_{\mathcal{H}}&\leq2 a(l(u_1)) |\nabla \bar{u}|^2
+2L_{a}(R)|l| |\bar{u}| |\nabla u_2||\nabla\bar{u}| \\[0.7ex]
&~~-2<f(u_1)-f(u_2),\bar{u}>_{L^{p,q}}-2(((\bar{\eta})^{\prime},\bar{\eta}%
))_{\mu} \\[1.0ex]
&\leq -2m|\nabla \bar{u}|^2+2L_{a}(R)|l||\bar{u}| |\nabla u_2||\nabla \bar{u}%
| \\[1.0ex]
&~~-2<f(u_1)-f(u_2),\bar{u}>_{L^{p,q}}-2(((\bar{\eta})^{\prime},\bar{\eta}%
))_{\mu} \\[1.0ex]
&\leq -2m|\nabla \bar{u}|^2+2m |\nabla \bar{u}|^2+\frac{1}{2m} L_a^2(R)|l|^2|%
\bar{u}|^2\|u_2\|^2 \\[1.0ex]
&~~-2<f(u_1)-f(u_2),\bar{u}>_{L^{p,q}}-2(((\bar{\eta})^{\prime},\bar{\eta}%
))_{\mu} \\[1.0ex]
&\leq \frac{1}{2m}L_a^2(R)|l|^2\|\bar{z}\|_{\mathcal{H}}^2\|u_2\|^2
-2<f(u_1)-f(u_2),\bar{u}>_{L^{p,q}}-2(((\bar{\eta})^{\prime},\bar{\eta}%
))_{\mu}, \\[1.0ex]
\end{split}%
\end{equation}
where $<\cdot,\cdot>_{L^{p,q}}$ is the duality between $L^{2p}$ and $L^q$.
The previous calculation is obtained formally taking the product in $%
\mathcal{H}$ between $\bar{z}$ and the difference of (\ref{eq3.5}) with $z_1$
and $z_2$ in place of $z$, and it can be made rigorous with the use of
mollifiers, see \cite[Theorem, Section 2]{C}. 
In fact, integrating by parts and by the fact that $\mu^{\prime} <0$ (see
again \cite[Section 2]{C}), we have 
\begin{equation*}
2(((\bar{\eta})^{\prime},\bar{\eta}))_{\mu}=-\lim_{s\rightarrow
0}\mu(s)|\nabla \bar{\eta}^t(s)|^2-\int_0^{\infty}\mu^{\prime}(s) |\nabla%
\bar{\eta}^t(s)|^2ds \geq 0.
\end{equation*}
Hence, the last term of the right hand side of (\ref{eq3.20}) can be
neglected.

At last, from (\ref{eq3.4}) we know that $f(u)$ is increasing for $|u|\geq M$%
, for some $M>0$. Fix $t\in (\tau,T]$, and let 
\begin{equation*}
\Omega_1:=\{x\in \Omega:|u_1(t,x)|\leq M~~\mbox{and}~~|u_2(t,x)|\leq M\},
\end{equation*}
and 
\begin{equation*}
N=2\sup_{|s|\leq M}|f^{\prime}(s)|.
\end{equation*}
Let $x\in\Omega_1$, then we have 
\begin{equation*}
2|f(u_1(x))-f(u_2(x))|\leq N|\bar{u}(x)|.
\end{equation*}
Then, by the monotonicity of $f(u)$ for $|u|\geq M $ and the Poincar\'e
inequality, it follows that 
\begin{equation}  \label{eq3.22}
\begin{split}
-2<f(u_1)-f(u_2),\bar{u}>_{L^{p,q}}&\leq
-2\int_{\Omega_1}(f(u_1(x))-f(u_2(x)))\bar{u}(x)dx \\
&\leq \int_{\Omega_1}N|\bar{u}(x)|^2dx \\[1.0ex]
&\leq N\|\bar{z}\|_{\mathcal{H}}^2.
\end{split}%
\end{equation}
(\ref{eq3.20})-(\ref{eq3.22}) imply that 
\begin{equation*}
\frac{d}{dt}\|\bar{z}\|^2_{\mathcal{H}}\leq \left( \frac{1}{2m}%
L_a^2|l|^2\|u_2\|^2+N\right)\|\bar{z}\|^2_{\mathcal{H}}.
\end{equation*}
The uniqueness and continuous dependence on initial data of solution to
problem \eqref{eq3.5} follow from the Gronwall inequality. Till now, we
finish the proof of the first assertion.

\textbf{Step 6.} (Further regularity) We are going to study further
regularity of $(u,\eta)$. {To this end, let us first consider the linear
operator $\mathcal{I}: L^2_{V\cap H^2(\Omega)}\rightarrow L^2_\mu(\mathbb{R}%
^+;D(V))$ defined by 
\begin{equation*}  \label{j1}
(\mathcal{I}\phi)(s)=\int_{-s}^0\phi(r)\,dr,\quad s\in\mathbb{R}^+.
\end{equation*}
Then, the operator $\mathcal{I}$ defined above is a linear and continuous
mapping. In particular, there exists a positive constant $K_{\mu}$, which is
the same as in Lemma \ref{lemma3.1}, such that, for any $\phi\in L_{V\cap
H^2(\Omega)}^2$, it holds 
\begin{equation*}  \label{Kmu1}
\|\mathcal{I}\phi\|^2_{L_{\mu}^2(\mathbb{R}^+; D(A))}\leq
K_\mu\|\phi\|^2_{L^2_{V\cap H^2(\Omega)}}.
\end{equation*}%
} 

Next, multiplying $(\ref{eq3.3})_1$ by $-\Delta u$ with respect to the inner
product of $H$, the Laplacian of $(\ref{eq3.3})_2$ by $\eta$ with respect to
the inner product of $L^2_{\mu}(\mathbb{R}^+;D(A))$, and adding the two
terms, we obtain 
\begin{equation}  \label{eq3.23}
\frac{d}{dt}\|z\|_{\mathcal{V}}^2+2a(l(u))|\Delta
u|^2+2(((\eta^t,(\eta^t)^{\prime})))_{\mu}=2(-f(u)+g,\Delta u).
\end{equation}
Since $f$ is a polynomial of odd degree, there exists a constant $d_0>0$,
such that 
\begin{equation}  \label{eqd}
f^{\prime}(u)\geq -\frac{d_0}{2},\qquad \forall u\in\mathbb{R}.
\end{equation}
Then, it follows from the above inequality, (\ref{eq3.4}), the Green formula
and the Young inequality that 
\begin{equation*}
\begin{split}
2(f(u),\Delta u)&=2\int_{\Omega}f_{2p-1}\Delta udx-2\int_{\Omega}
f^{\prime}(u)\nabla u\cdot \nabla udx \\
&\leq \frac{2}{m}f_{2p-1}^2|\Omega|+\frac{m}{2}|\Delta u|^2+d_0|\nabla u|^2.
\end{split}%
\end{equation*}
Again by the Young inequality, we have 
\begin{equation*}
2(g,\Delta u)\leq \frac{m}{2}|\Delta u|^2+\frac{2}{m}|g|^2.
\end{equation*}
Together with \eqref{eq2.2}, (\ref{eq3.23}) becomes 
\begin{equation}  \label{eq3.24}
\frac{d}{dt}\|z\|^2_{\mathcal{V}}+m|\Delta
u|^2+2(((\eta^t,(\eta^t)^{\prime})))_{\mu}\leq \Theta,
\end{equation}
where we have used the notation $\Theta=\frac{2}{m}f_{2p-1}^2|\Omega|+d_0|%
\nabla u|^2+\frac{2}{m}|g|^2$, which belongs to {$L^1(\tau,T) $.} Under the
suitable spatial regularity assumptions on $\eta$, integration by parts in
time and using $(h_1)$, we obtain 
\begin{equation*}
(((\eta^t,(\eta^t)^{\prime})))_{\mu}=-\int_0^{\infty}\mu^{\prime}(s)|\Delta
\eta^t(s)|^2ds\geq 0.
\end{equation*}
Therefore, the term $2(((\eta^t,(\eta^t)^{\prime})))_{\mu}$ in (\ref{eq3.24}%
) can be neglected, we integrate (\ref{eq3.24}) between $\tau$ and $t$,
where $t \in(\tau,T)$, which leads to 
\begin{equation}  \label{eq3.25}
\|z(t)\|^2_{\mathcal{V}}+m\int_{\tau}^t |\Delta u(s)|^2ds\leq \|z(\tau)\|^2_{%
\mathcal{V}}+ \int_{\tau}^t\Theta(s)ds.
\end{equation}
From the above estimation, we conclude that 
\begin{equation*}
\begin{split}
& u\in L^{\infty}(\tau,T,V)\cap L^2(\tau,T;D(A)); \\
& \eta\in L^{\infty}(\tau,T;L^2_{\mu}(\mathbb{R}^+;D(A))).
\end{split}%
\end{equation*}
Concerning the assertion $(ii)$ of this theorem, the continuity of $u$
follows again using a slightly modified version of \cite[Lemma III.1.2]{T1}.
The continuity of $\eta$ can be proved mimicking the idea of the proof of
Step 4 of $(i)$, with $D(A)$ in place of $V$. The proof of this theorem is
complete. $\Box$

\end{document}